\documentclass[onefignum,onetabnum]{siamart190516}

\usepackage{lipsum}
\usepackage{amsfonts}
\usepackage{graphicx}
\usepackage{epstopdf}
\usepackage{algorithmic}
\ifpdf
  \DeclareGraphicsExtensions{.eps,.pdf,.png,.jpg}
\else
  \DeclareGraphicsExtensions{.eps}
\fi

\newsiamremark{remark}{Remark}
\newsiamremark{hypothesis}{Hypothesis}
\crefname{hypothesis}{Hypothesis}{Hypotheses}
\newsiamthm{claim}{Claim}

\headers{Dynamic Optimization with Convergence Guarantees}{M. P. Neuenhofen and E. C. Kerrigan}

\title{Dynamic Optimization with Convergence Guarantees}

\author{Martin P. Neuenhofen\thanks{Department of Electrical \& Electronic Engineering, Imperial College London, SW7~2AZ London, UK,
  (\email{m.neuenhofen19@imperial.ac.uk}, \url{http://www.MartinNeuenhofen.de}).}
\and Eric C. Kerrigan\thanks{Department of Electrical \& Electronic Engineering and Department of Aeronautics, Imperial College London, SW7~2AZ London, UK, 
  (\email{e.kerrigan@imperial.ac.uk}, \url{http://www.imperial.ac.uk/people/e.kerrigan}).}}

\usepackage{amsopn}

\ifpdf
\hypersetup{
	pdftitle={Direct Transcription of Dynamic Optimization Problems using Integral Penalties and Barriers},
	pdfauthor={M. N. Neuenhofen and E. C. Kerrigan}
}
\fi

\usepackage{amsmath, amsfonts}
\usepackage{color}
\usepackage{amssymb}
\usepackage{mathtools}
\usepackage{csquotes}

\usepackage{calc}

\theoremstyle{plain}
\newtheorem{thm}{Theorem}

\newtheorem{lem}{Lemma}
\newtheorem{prop}{Proposition}
\theoremstyle{plain}
\newtheorem{defn}{Definition}

\theoremstyle{plain}
\newtheorem{rem}{Remark}

\usepackage{mdframed}

\usepackage{todonotes}

\usepackage{cite}
\usepackage{comment}

\usepackage{graphicx}
\usepackage{subfigure}          

\usepackage{placeins}

\usepackage{tabularx}

\usepackage{adjustbox}

\usepackage{lmodern,textcomp}

\usepackage{enumerate}

\usepackage{xspace}

\usepackage{microtype}

	\newcommand{\inv}{^{-1}\xspace}

	\newcommand{\blambda}{\boldsymbol{\lambda}\xspace}

	\newcommand{\bs}{\textbf{s}\xspace}
	\renewcommand{\t}{^\textsf{T}\xspace}

	\newcommand{\away}[1]{}

	\newcommand{\nx}{n_x}
	\newcommand{\nquad}{{n_{\text{q}}}\xspace}
	\newcommand{\nX}{{N_{X}}\xspace}

	\newcommand{\nc}{{n_c}\xspace}
	\newcommand{\nb}{{n_b}\xspace}

	\newcommand{\R}{\mathbb{R}\xspace}

	\newcommand{\N}{\mathbb{N}\xspace}

	\newcommand{\cV}{\mathcal{V}\xspace}

	\newcommand{\cB}{\mathcal{B}\xspace}

	\newcommand{\cL}{\mathcal{L}\xspace}
	\newcommand{\cX}{\mathcal{X}\xspace}
	\newcommand{\cP}{\mathcal{P}\xspace}

	\newcommand{\cU}{\mathcal{U}\xspace}
	
	\newcommand{\cC}{\mathcal{C}\xspace}

	\newcommand{\cT}{\mathcal{T}\xspace}
	\newcommand{\cO}{\mathcal{O}\xspace}

	\newcommand{\bD}{\textbf{D}\xspace}
	\newcommand{\bS}{\textbf{S}\xspace}

	\newcommand{\bx}{\textbf{x}\xspace}

	\newcommand{\by}{\textbf{y}\xspace}

	\newcommand{\be}{\mathbf{1}\xspace}
	\newcommand{\bei}[1]{{\textbf{e}}\xspace}

	\newcommand{\bO}{\textbf{0}\xspace}

	\newcommand{\opdiag}{\operatorname{diag}\xspace}

	\newcommand{\tol}{{\textsf{tol}}\xspace}

	\newcommand{\commentout}[1]{}
	
	\newcommand{\infexpr}{{\operatornamewithlimits{min}_{x_h \in \cX_{h,p}}\|x^\star_{\omega,\tau}-x_h\|_\cX}\xspace}

	\newcommand\myeqLH{\mathrel{\stackrel{\makebox[0pt]{\mbox{\normalfont\tiny L'H}}}{=}}}

\usepackage{datetime,url,balance}

\usepackage{soul}
\usepackage{color}

\DeclareRobustCommand{\submithla}[1]{{\textcolor{black}{#1}}}

\newlength\colA
\newlength\colB
\newlength\colC
\newlength\colD
\newlength\colE
\begin{document}
	
	
	
	\maketitle
	
	\begin{abstract}
		We present a novel direct transcription method to solve optimization problems subject to nonlinear differential and inequality constraints.
		We prove convergence of our numerical method under reasonably mild assumptions: boundedness and Lipschitz-continuity of the problem-defining functions.
		We do not require uniqueness, differentiability or constraint qualifications to hold and we avoid the use of Lagrange multipliers. Our approach differs fundamentally from well-known methods based on collocation;
		we follow a penalty-barrier approach, where we compute integral quadratic penalties on the equality path constraints and point constraints, and integral log-barriers on the inequality path constraints. The resulting penalty-barrier functional can be minimized numerically using finite elements and penalty-barrier interior-point nonlinear programming solvers.
		Order of convergence results are derived, even if components of the solution are discontinuous.
		We also present numerical results to compare our method against collocation methods. The numerical results show that for the same degree and mesh, the computational cost is similar, but that the new method can achieve a smaller error and  converges in cases where collocation methods fail.
	\end{abstract}
	
	
	
	\section{Introduction}
	
	\subsection{An Important Class of Dynamic Optimization Problems}
	\label{sec:classDOP}
	Many optimal control, estimation, system identification and design problems can be written as a dynamic optimization problem in the Lagrange form
	\newcommand{\refOCP}{(DOP)\xspace}
\begin{align}
 	\operatornamewithlimits{min}_{x:=(y,z) \in \cX} 	
    \int_\Omega f( \dot{y}(t),y(t),z(t),t)&\, \mathrm{d}t  \label{eqn:OCP} \tag{DOPa}\\
	\text{subject to\quad}  	
    b\left(y(t_1),\ldots,y(t_M)\right) &=0,\label{eqn:point} \tag{DOPb}\\
      c\left(\dot{y}(t),y(t),z(t),t\right)
      &=0\quad\text{ f.a.e.\ }t \in \Omega, \label{eqn:dae} \tag{DOPc}\\
      z(t) &\geq 0 \quad\text{ f.a.e.\ }t \in \Omega,\label{eqn:positive} \tag{DOPd}
\end{align}
	where the open bounded interval
	$ 	\Omega := (t_0,t_E) \subsetneq \R$, $\cX$ is an appropriately-defined Hilbert space for solution candidates $x:=(y,z)$ such that $y$ is continuous, and \enquote{f.a.e.} 
	means \enquote{for almost every} in the Lebesgue sense. Detailed definitions and assumptions are given  in Section~\ref{sec:assumptions}. An optimal solution will be denoted with $x^\star$.
	We note that the form~\refOCP is quite general and adopted here to minimize notation.

	Ordinary differential equations (ODEs) and path constraints are included via the differential-algebraic equations (DAE) in~\eqref{eqn:dae} and the inequalities~\eqref{eqn:positive}. The point constraints~\eqref{eqn:point} enforce boundary constraints, such as initial or final values on the state $y$, or include values obtained by measurements at given time instances.

	With techniques presented in \cite[Sect.~3.1]{Patterson:2014:GMS:2684421.2558904},\cite[Sect.~2--3]{BrysonWildly},\cite[Chap.~4]{Betts2nd},
	problems in the popular Bolza or Mayer forms with general inequalities and free initial- or end-time can be converted into the form~\refOCP; see Appendix~\ref{sec:convert}. In turn, many problems from control,  estimation and system identification can be stated in Bolza or Mayer form~\cite{Betts2nd}.

	
	
	Problem~\refOCP is infinite-dimensional, because the optimization is over function spaces subject to an uncountable set of constraints. It is  very hard or impossible to compute an analytic solution, in general. Hence, one often has to resort to numerical methods to solve \refOCP.
	When doing so, it is important to eliminate whether  features of the numerical solution have arisen from physical principles or numerical failure.
	The need for a numerical method, which has a rigorous proof that the numerical solution convergences to the exact solution, is therefore essential in practice.
	
	One of the more established choices for the numerical solution of \refOCP is to discretize via
	direct collocation and finite elements~\cite{Conway2012,Betts2nd,rao_asurvey,Patterson,KellyMatthew}. 
	Recall that explicit Runge-Kutta methods are unsuitable for stiff problems and that many popular implicit methods for solving differential equations, e.g.\ variants of Gauss schemes, can be interpreted as collocation methods. Collocation methods include certain classes of implicit Runge-Kutta, pseudospectral, as well as Adams and backward differentiation formula methods~\cite{KellyMatthew,Betts2nd,rao_asurvey,arevelo:2002}.
	However, as is known \cite[Sect.~2.5\,\&\,4.14]{Betts2nd},\cite{ChenBiegler16,kameswaran_biegler_2008}, collocation methods can fail to converge if care is not taken. In \cite{CDC20PBF} we present an example where three commonly used collocation-based direct transcription methods diverge, and below in Section~\ref{sec:numExp:CQP} we give a parametric problem for which Legendre-Gauss-Radau collocation~\cite{Betts2nd} of any degree rings.
	
	Notice that \refOCP includes problems with mixed differential and inequality path constraints, for which indirect methods~\cite{Conway2012,Betts2nd} have only a limited range of applicability. Even when applicable, indirect methods require sophisticated user knowledge to set up suitable co-state and switching-structure estimates \cite{Boehme2017}. A detailed discussion of available methods in the literature is given in Section~\ref{sec:literatureReview}.
	
	There is a scarcity of rigorous proofs that show that high-order collocation schemes for dynamic optimization methods converge to a feasible or optimal solution as the discretization is refined. The assumptions in the literature are often highly technical, difficult to enforce or not very general.
	
	\subsection{Contributions}
	
    This paper aims to overcome the limitations of the  numerical methods mentioned above by presenting a novel direct transcription method  for solving~\refOCP. 
    Our method combines the following ingredients: quadratic integral penalties for the equality (path) constraints; logarithmic integral barriers for the inequality path constraints; and direct numerical solution via finite elements.
	It is this combination, together with a rigorous proof of convergence, that amounts to a novel direct transcription method. 
	We also provide order-of-convergence results.
    
	As detailed in Section~\ref{sec:assumptions}, we  only require existence of a solution to~\refOCP and mild assumptions on the boundedness and Lipschitz continuity of $f,c,b$.
	In contrast to existing convergence results:
	\begin{itemize}
		\item The solution $x^\star$ 
		does not need to be  unique.
		\item  $f,c,b$ can be non-differentiable everywhere.
		\item  We do not require the satisfaction of a constraint qualification for the discretized finite-dimensional optimization problem, such as the Linear Independence Constraint Qualification (LICQ), Mangasarian-Fromovitz Constraint Qualification (MFCQ) or Second-Order Sufficient Conditions (SOSC).
		\item Uniqueness or global smoothness of states or co-states/ adjoints do not need to hold.
		\item Local uniqueness assumptions, as in~\cite{arXiv:2017}, are removed.
	\end{itemize}

	The idea behind our new, Penalty-Barrier-Finite Element method (PBF), is to minimize the following unconstrained penalty-barrier function
	\begin{align}
	    \Phi(x) := F(x) + \frac{1}{2\cdot \omega} \cdot r(x) + \tau \cdot \Gamma(x), \label{eqn:def:Phi}
	\end{align} 
	where
	\begin{align}
	F(x) := &\int_\Omega \, f\left(\dot{y}(t),y(t),z(t),t\right)\mathrm{d}t \label{eqn:def:F}\end{align}
	is the objective,
	\begin{align}
		\begin{split}
		r(x) := &\int_\Omega\|c\left(\dot{y}(t),y(t),z(t),t\right)\|_2^2\,\mathrm{d}t + \|b\left(y(t_1),y(t_2),\ldots,y(t_M)\right)\|_2^2
		\end{split}\label{eqn:def:r}
	\end{align}
	is the integral quadratic penalty for the equality path- and point constraints, and
	\begin{align}
		\Gamma(x):=&-\sum_{j=1}^{n_z} \int_\Omega \log\big( z_{[j]}(t)\big)\,\mathrm{d}t\label{eqn:def:Gamma}
	\end{align}
	is an integral logarithmic barrier for the inequality path constraints. We provide an analysis that shows that one can construct trajectories $x_h$ that converge in the following \textit{tolerance-accurate} sense: the \textit{optimality gap}
	\begin{align}
	    g_\text{opt} :=& \max\lbrace 0,F(x_h)-F(x^\star)\rbrace
	\end{align}
	and \textit{feasibility residual}
	\begin{align}
	    r_\text{feas} :=& r(x_h)
    \end{align}
	converge to zero as the discretization mesh becomes finer and the parameters $\tau,\omega>0$ converge to zero. Order-of-convergence results will specify the rate at which $g_\text{opt}$ and $r_\text{feas}$ approach zero.
	
	The above functions \eqref{eqn:def:F}--\eqref{eqn:def:Gamma} look similar to those encountered in well-known finite-dimensional penalty-barrier methods. However, in order to deal with the infinite-dimensional nature of the problem, note the use of \emph{integrals} in the penalty \emph{and} barrier terms. If the problem had been finite-dimensional in $x$ and if $r$ had been the squared $2$-norm of finitely many equality constraints,  then it would be given that the minimizer of $\Phi$  converges to the solution $x^\star$ under mild assumptions as  $\tau,\omega$ converge to zero~\cite{SUMT}. The infinite-dimensional case considered here, however, is more involved and requires a careful analysis relating $\tau,\omega$ to  parameters of the  discretization. This is because once we discretize on a mesh  and seek to compute an approximate solution $x_h$ on the mesh, the degrees of freedom for $x_h$ depend on the size of the finite element space. If we were to draw an analogy with the finite dimensional case, then the equivalent number of equality constraints depends on the number of quadrature points for numerically evaluating the integral in $r$. If $\omega$ is too large then~$x_h$ will not converge to satisfying the equivalent set of equality constraints. If $\omega$ is too small with respect to the mesh size, then there are not enough degrees of freedom, resulting in a potentially feasible but suboptimal  solution~\cite[p.~1078]{Hager90}.
	The effects of $\tau$ are more intricate, since they relate to a local Lipschitz property of $\Phi$ that is relevant for the stability of the finite element discretization.
	A balance must thus be taken between the size of the finite element space, the quadrature rule and the parameters $\omega,\tau$. This requires a non-trivial analysis, which is the core contribution of the paper. 
	
	For the special case when $f=0$ and the inequality  constraints~\eqref{eqn:positive} are removed, \refOCP reduces to finding a solution of the differential equations~\eqref{eqn:dae} subject to~\eqref{eqn:point}. Our method then reduces to that of the the least-squares method~\cite{locker1978,ascher1978}, i.e.\ solving a sequence of finite-dimensional least-squares problems where the  mesh is refined until the feasibility residual $r_\text{feas}$ is below the required tolerance. 
	The method in this paper can therefore be interpreted as a generalization of least-squares methods for differential equations to   dynamic optimization problems. The extension and analysis is not straightforward, because of the interaction in $\Phi$ between $r$ and the additional two terms involving the objective and inequality constraints.

	\subsection{Motivation from Collocation Methods}
We motivate our method from the perspective of collocation methods.

A desirable method for solving optimal control problems is Legendre-Gauss-Radau collocation because it is easy to implement (and high-order consistent): the method constructs piecewise polynomials (of high degree) using a nodal basis and solves the path constraints at a finite number of points. The nodal basis values are determined by solving a large sparse nonlinear program.

However, for solutions with singular arcs, which occur in a large number of applications, the numerical solutions can ``ring''~\cite[Sect.~4.14.1]{Betts2nd}. In particular, the solution polynomial and the residuals of the path constraints will oscillate between the collocation points --- that is, the path residuals will not converge to zero everywhere. A remedy is regularization: a small convex quadratic term is added to the objective to penalize numerical noise. Unfortunately, in most cases this remedy does not work because either the penalty is too small to remove all noise or so large that it alters the problem's solution.

The idea with the penalty method is to make ringing impossible by adding collocation points inbetween the original collocation points, where otherwise the states, controls and residuals could ring. The theoretical vehicle for this approach are integrals and penalties. Integrals, once discretized by means of numerical quadrature, can be expressed with a set of weights and abscissae, alias collocation points. Penalties, in replacement for exact constraints, will prevent any issues related to the ratio between the number of degrees of freedom and the number of constraints, such as over-determination. The resulting scheme remains easy to implement while effectively forcing global convergence of the path constraints --- as we rigorously prove in the remainder of the paper. In particular, we prove that the feasibility residual converges to zero.

We stress that the integral penalty and log-barrier provide a useful natural scaling for the NLP. This is certainly desirable from a computational perspective, because experience shows that numerical treatment of an NLP depends significantly on scaling~\cite[Chap.~1.16.5,~4.8]{Betts2nd}. The large-scale methods in \cite{ForsgrenGill,ALM_IPM} use a merit function that treats equality constraints with a quadratic penalty and inequality constraints with a log-barrier term. Typically, as the NLP becomes larger, caused by a finer discretization, the NLP becomes more challenging to solve, in the sense that the number of iterations to converge increases. In contrast, in this paper the NLP merit function matches  the infinite-dimensional merit function in the limit, which mitigates numerical issues that might otherwise arise.

\commentout{
We stress why the integral penalty and log-barrier provide a useful natural scaling for the NLP: It makes sense to solve the NLP with a (primal-dual) interior-point method, such as IPOPT. The iterate of the solver can be associated with the minimizer of a penalty-barrier merit function. Consider for instance the method in \cite{ForsgrenGill,ALM_IPM}, that treat equality constraints with a quadratic penalty and inequality constraints with a log-barrier. Typically, as the NLP becomes larger, caused by a finer discretization, it becomes more challenging to solve, in the sense that the number of iterations to converge increase.

Here instead, by careful formulation of the infinite-dimensional problem with integrals, we have the opportunity to identify the NLP merit function with a convergent discretization of an infinite-dimensional penalty-barrier functional of indeed the optimal control problem that we wish to solve. If we were to minimize that penalty-barrier functional in the SUMT approach with a primal method --which could be called an infinite-dimensional interior-point method-- then clearly the iteration sequence is mesh-independent, as there is no mesh. Hypothetically, for a sufficiently fine discretization of the penalty-barrier functional, each individual finite-dimensional iterate of a primal method on the discretized penalty-barrier functional should converge to its respective infinite-dimensional iterate of the infinite-dimensional primal method. This would imply convergence of the NLP iterates at a mesh-independent rate. Less hypothetical, the integral formulation on the infinite dimension induces a natural scaling for the (finite-dimensional) NLP.
This is certainly desirable from a computational perspective because experience shows that numerical treatment of NLP depends significantly on scaling~\cite[1.16.5,~4.8]{Betts2nd}.}

\subsection{Difference to Galerkin Methods for Optimal Control}
There can be confusion on the relation between integral penalty methods and Galerkin-type finite element methods. Galerkin-type methods solve operator equations of the form $T(u)=0$, where $u$ is a function and $T$ is a functional. Defining a solution space $\cU$ and a test space $\cV$, Galerkin methods seek solutions $u$ that satisfy
\begin{align}
	\text{Find }u \in \mathcal{U}\text{ such that } \langle v,T(u)\rangle=0\quad\forall v \in \cV\,, \label{eqn:Galerkin}
\end{align}
for a scalar product over a super-space of $\cU,\cV$. Approximations to $u$ can be computed in finite time by approximating $\cU,\cV$ with finite-dimensional spaces, e.g.\ finite element spaces. There are three important kinds of methods: Galerkin, where $\mathcal{V}=\mathcal{U}$; Galerkin Least-Squares, where $\cV=T(\cU)$; and collocation, where $\cV$ is a set of dirac delta functions.

If $T$ is affine, then the scalar product can be rephrased as
\begin{align*}
	\text{Find }u \in \mathcal{U}\text{ such that } a(u,v)=F(v)\ \forall v \in \mathcal{V}\,.
\end{align*}
For $a$ bi-linear and $F$ linear, the equation takes on the form $v\t A u =v\t f$, such that upon a finite basis for $\mathcal{U,V}$ a linear equation system results for the coefficient vector of $u$ in that basis. Sufficient conditions for $u$ to be unique exist (e.g.\ coercivity).

Galerkin, i.e.\ $\cV=\cU$, is the most natural choice when $T$ is the first variation of a functional. To see this, consider the convex objective $J(u)=0.5\,a(u,u)-F(u)$, whose minimizer $u^\star$ satisfies $\langle v,T(u)\rangle \equiv a(u^\star,v)-F(v)=0 \ \forall v \in \cU$ when $T$ is the functional derivative of $J$, i.e.\ $T=\delta J$, hence $\langle T(u),v \rangle =\delta J(u;v)$. Galerkin Least-Squares results in the normal equations, forcing symmetry while not necessarily coercivity (e.g.\ when $A$ is singular). In the nonlinear case, $A$ is the Jacobian of $T$.
Collocation adds no natural benefit other than that it is easy to implement, particularly when $T$ is nonlinear, because then $A$ and $A \cU$ would constantly change, making Galerkin and Galerkin Least-Squares otherwise intractably expensive for implementation.

Galerkin methods are tailored for solving linear or mildy nonlinear operator equations (e.g.\ Maxwell, Navier-Stokes, Euler). Collocation methods are preferred in optimal control because the setting is entirely different to solving PDEs: In optimal control one minimizes highly nonlinear functionals instead of solving mildly nonlinear operator equations. Nevertheless, defining $J=\Phi$ and $T=\delta \Phi$, since we are minimizing $\Phi$, the optimal control problem could in principle be expressed as an operator problem and solved via Galerkin-type methods. Since one aims to minimize $J$, the necessary optimality condition $\langle T(u),v \rangle =\delta J(u;v)=0 \ \forall v\in\cU$ must be satisfied, as is solved by Galerkin. However, in practice, for reasons of efficiency and robustness, the minimizer $u$ of $J$ is rather constructed not via Galerkin \eqref{eqn:Galerkin} (possibly solved via Newton iterations), but with numerical quadrature and powerful modern nonlinear optimization software packages, such as IPOPT \cite{IPOPT}. The former method is indirect, the latter is direct. The latter is preferable due to the  disadvantages of the indirect approach for highly nonlinear problems, as discussed at the end of Section~\ref{sec:Lit2}.

We discuss the least-squares character of $\Phi$ and contrast this against Galerkin Least-Squares methods. Suppose $b=0$. If $\Phi(x)=r(x)$, then the necessary optimality conditions of minimizing $\Phi$ are equivalent to the Galerkin Least-Squares equations for solving $T(x)=0$ for $T(x):=c(\dot{y},y,z,t)$. Choosing finite elements $\cU=\cX_h$, defined in Section~\ref{sec:FEM}, one obtains $n:=\operatorname{dim}(\cX_h)$ equations for $n$ unknown finite element coefficients of $x_h \in \cX_h$. We give an example in \cite{CDC20MALM} of an optimal control problem where $x_h$ is uniquely determined when minimizing $r(x)$. Hence, utilization of Galerkin Least-Squares for constraint treatment would neglect optimality w.r.t.\ the objective $F$.

In summary, we present here a novel method that cannot be classified as a variant of Galerkin Least Squares or Galerkin Collocation. Our method is only a Galerkin method in the same weak sense as any FEM in optimal control is, since they all discretize states and controls by means of a finite element space and minimize functionals over that finite element space.

\subsection{Literature on Integral Penalty Transcription Methods}
\label{sec:literatureReview}
	
	Quadratic integral penalties were introduced in~\cite{Courant43} within a Rayleigh-Ritz method for partial differential equations with homogeneous Dirichlet boundary conditions by quadratic penalization. In \cite{Russell65}, existence and convergence of solutions to two integral penalty functions of generic form were analyzed, but without a discretization scheme.
	
	Quadratic penalty functions of a less generic form, suiting optimal control problems with explicit initial conditions, were studied in~\cite{Balakrishnan68}. The analysis focuses on the maximum principles that arise from the penalty function and their connection (under suitable assumptions on smoothness and uniqueness) to Pontryagin's maximum principle, laying groundwork for an indirect solution to the penalty function. A numerical method is not proposed. In~\cite{Jones70}, the analysis is extended to inequality path constraints with a fractional barrier function. Using suitable smoothness and boundedness assumptions on the problem-defining functions, it is shown that the unconstrained minimizer converges from the interior to the original solution. As in~\cite{Balakrishnan68}, the analysis uses first-order necessary conditions. A numerical scheme on how to minimize the penalty functional is not presented. Limitations are in the smoothness assumptions.
	
	In \cite{DeJulio70} the penalty function of~\cite{Balakrishnan68} is used for problems with explicit initial conditions, linear dynamics and no path constraints. Under a local uniqueness assumption, convergence is proven for a direct discretization with piecewise constant finite elements. The approach is extended in~\cite{Hager90} to augmented Lagrangian methods with piecewise linear elements for the states. The analysis is mainly for linear-quadratic optimal control, which is used for approximately solving the inner iterations. The error in the outer iteration (augmented Lagrangian updates) contracts if the initial guess is sufficiently accurate~\cite[Lem.~3]{Hager90}. It is unclear in~\cite{Hager90} whether the proposal is to solve the inner iteration via direct or indirect discretization and, in case of the former, how the integrals are evaluated. Limitations are with the existence of a solution for each inner iteration, no path constraints for the states, sparsity of algorithmic details, and the construction of a sufficiently accurate initial guess for the outer iteration to converge.

\commentout{	
\subsection{Literature on Integral Penalty Methods}
\label{sec:literatureReview}
	
	Quadratic integral penalties were introduced as early as~\cite{Courant43} in the context of a Rayleigh-Ritz method, where they were used to solve partial differential equations with homogeneous Dirichlet boundary conditions by quadratic penalization; this can be considered an indirect finite element method. In \cite{Russell65} the existence and convergence of solutions to two integral penalty functions of generic form were analyzed, but without a discretization scheme.
	
	Quadratic penalty functions of a less generic form, which suit optimal control problems with explicit initial conditions, were studied in~\cite{Balakrishnan68}. The analysis focuses on the maximum principles that arise from the penalty function and their connection (under suitable assumptions on smoothness and uniqueness) to Pontryagin's maximum principle, laying groundwork for an indirect solution to the penalty function; a numerical method is not proposed. In~\cite{Jones70}  the analysis is extended for inequality path constraints by a fractional barrier function. Under certain smoothness and boundedness assumptions on the problem-defining functions it is shown that the unconstrained minimizer converges from the interior to the original solution. As in~\cite{Balakrishnan68}, the analysis uses first-order necessary conditions. A numerical scheme on how to minimize the penalty functional is not presented. Limitations are in the smoothness assumptions and lack of convergence-analysis for a numerical scheme to solve the function.
	
	In \cite{DeJulio70} the penalty function of~\cite{Balakrishnan68} is used for explicit initial conditions, linear dynamics and no path constraints. Under a local strong uniqueness assumption,  convergence  is proven for a direct discretization with piecewise constant finite elements. This approach is extended in~\cite{Hager90}   to an augmented Lagrangian method with piecewise linear elements for the states. As in~\cite{DeJulio70}, the analysis is mainly for linear-quadratic optimal control, which is used for approximately solving the inner iterations.  The error in the outer iteration (augmented Lagrangian update) contracts if the initial guess is sufficiently accurate~\cite[Lem.~3]{Hager90}. It is unclear in~\cite{Hager90} whether the proposal is to solve the inner iteration via direct or indirect discretization and, in case of the former, how the integrals are evaluated. Limitations are with the existence of a solution for each inner iteration, no path constraints for the states, sparsity of algorithmic details and the construction of a sufficiently accurate initial guess.}

\commentout{
	\subsection{Literature on Modern Direct Transcription Methods}\label{sec:Lit2}
	
	The literature on convergence of numerical methods for the solution of ODEs and DAEs is substantial. This is not the case for dynamic optimization.
	
	A convergence proof for a discretization based on the explicit Euler method is given in \cite{Maurer}, which makes the following strong assumptions: (i)  problem-defining functions must be locally differentiable with Lipschitz continuous derivatives; (ii) there must be a local solution where the trajectories of the state and free variables are continuously differentiable and continuous, respectively; (iii) a  homogeneity condition on active constraints; (iv) surjectivity of linearized equality constraints; and (v) a  coercivity assumption.
	These conditions are sophisticated, difficult to understand, and  hard to verify and ensure  by construction.
	The conditions assert that the first order optimality conditions of the infinite-dimensional problem result in a unique and convergent solution for the optimality conditions of the collocation method. This is why convergence proofs for other collocation schemes make similar assumptions.
	Furthermore, since Euler's method converges only of first order, convergence results based on Euler's method are of limited practical use compared to results for higher-order methods, if such results are available.

	In \cite{Kameswaran} the authors propose an $\ell_1$-penalty method for a high-order collocation method and demonstrate the stabilization effect of this approach for the numerical test problem~ of\cite{AlyChan} in Mayer form with path constraints. The result is experimental. The discussion in~\cite{Kameswaran} assumes similar conditions as in \cite{kameswaran_biegler_2008}, where they present a convergence result for problems without path constraints and do not use the $\ell_1$-penalty method from~\cite{Kameswaran}.
	The result in~\cite{kameswaran_biegler_2008} relies on a number of strong assumptions: (i) functions defining the problem must be sufficiently smooth; (ii) the state and co-state trajectories must be sufficiently smooth; (iii) the nonlinear program arising from the discretization must satisfy LICQ and SOSC.
	
	Another proof of high-order convergence for a direct collocation method is given in~\cite{RaoHager}. This proof considers the Mayer form with constraints. The authors show convergence of their scheme under the following assumptions:
	(i) the solution is locally unique;
	(ii) the states have a strong first derivative; (iii) the state and co-state trajectories  have two square integrable derivatives;
	(iv) the Hamiltonian satisfies a local strong convexity property;
	(v) the objective and the ODE function are twice Lipschitz differentiable.
	
	A convergence result for  a pseudospectral method for the control of constrained feedback linearizable systems is given in~\cite{gongEtAl:2006}. In order to provide a proof that does not require dualization, the strong assumption is made that the derivatives of the interpolating polynomial for the state converges uniformly to a continuous function, which implies that the optimal input has to be continuous. This assumption was relaxed in~\cite{kangEtAl:2007} to allow for discontinuous optimal inputs, under the assumptions that (i)~both the optimal state and input trajectories are  piecewise differentiable, and (ii)~the path constraints define a convex set.
	
	The assumption on feedback linearizable systems in~\cite{gongEtAl:2006,kangEtAl:2007} was relaxed in~\cite{gongEtAl:2008} to allow for more general nonlinear ordinary differential equations. However, the following strong assumptions are made: (i)~the problem-defining functions are continuously differentiable, (ii) the  gradients of the functions are Lipschitz continuous, (iii)~the optimal state trajectory is continuously differentiable, which requires that the optimal control trajectory be continuous.
	
	Direct multiple shooting~\cite{rao_asurvey,Betts2nd,SOSC} is an  effective method for solving dynamic optimization problems.  Convergence proofs for  shooting methods  for quadratic regulator problems with linear dynamics and linear inequality constraints are presented in~\cite{pannochiaEtAl:2010,pannochiaEtAl:2015,yuzEtAl:2005}, where the  differential equation  is solved using matrix exponentials.
	
	A high-order convergence proof for explicit, fixed-step size Runge-Kutta methods is given in~\cite{SchwartzPolak:1996}, which  extends  some of the results presented in~\cite{PolakBook1997} for the Euler method. The assumptions in~\cite{SchwartzPolak:1996} are: (i)~the inputs are constrained to lie in a convex and compact norm-ball, (ii)~the state derivatives can be written as an explicit,  continuously differentiable function of the state and input, and (iii)~that the differential equations are time-invariant.
	
	Convergence results for Runge-Kutta methods for unconstrained optimal control problems are available in \cite{Hager2000}, subject to the following assumptions: (i) the optimal states and their first and second derivatives are globally bounded, (ii) the optimal control is continuously differentiable, (iii) functions defining the problem are twice differentiable, and (iv) a local coercivity property of the linear-quadratic approximation in the local minimizer.
	
	Indirect methods have also been widely studied for computing the solution of ~\refOCP~\cite{rao_asurvey,Betts2nd}. In these, the calculation of variations is used to determine the optimality conditions for the optimal arcs.
	Indirect methods have been less successful in practice than direct methods.
	Firstly, the resulting optimality system, namely a Hamiltonian boundary-value problem (HBVP), is difficult to solve numerically. There are robustness issues when no accurate initial guess for the solution of the HBVP is given \cite{BrysonWildly}. Secondly, for complicated problems it is difficult or impossible to determine the optimality conditions to begin with, hence there are problems which cannot be solved by indirect methods \cite{Boehme2017}. Even when possible, state-constrained problems require estimates of the optimal switching structure. Workarounds for state constraints, such as by saturation functions, have been proposed in~\cite{WANG20174070}. Thirdly, for singular-arc problems, optimality conditions of higher order need to be used in order to determine a solution \cite[Sec.~1.4]{LAMNABHILAGARRIGUE1987173}. Finally, the optimality conditions can be ill-posed, for example when the co-state solution is non-unique. This always happens when path constraints are linearly dependent, because a dual solution is then non-unique.}
	
	\subsection{Literature on Non-Penalty Transcription Methods}\label{sec:Lit2}
	\subsubsection{Direct Methods}
	
	Compared to the solution of ODEs and DAEs, the literature on convergence of numerical methods for dynamic optimization is relatively sparse.  Table~\ref{tab:literatureReview} contains an overview of some of the convergence results, which are discussed in more detail in the following.
	
\begin{table}[htpb!]
\label{tab:literatureReview}
\caption{Assumptions for convergence results in the literature. Cells are empty when no assumptions in this category are made or when assumptions in this category cannot be matched to our problem formulation. Abbreviations: conv.=convex; cont.=continuous; suff.=sufficiently; hom. cond.=homogeneity condition; surj. lin. eq. constr.=sujective linear equality constraints; RKM=Runge-Kutta method}
\settowidth\colA{$fc$ L-cont. in $y$}
\settowidth\colB{$\nabla\lbrace f,c,b \rbrace$}
\settowidth\colC{explicit ODEx}
\settowidth\colD{LICQ, SOSC}
\settowidth\colE{no path constraints}
\tabcolsep=4pt
\centering
\begin{adjustbox}{angle=90}
\begin{tabular}{@{}|l||>{\raggedright}p{\colA}|>{\raggedright}p{\colB}|>{\raggedright}p{\colC}|>{\raggedright}p{\colC}
|>{\raggedright}p{\colD}|>{\raggedright}p{\colE}||p{\colE}|@{}}
\hline
Ref. 					& $f,c,b$    		& $\nabla\lbrace f,c,b \rbrace$ 						& $y$                             	& $z$                          		& NLP        				& remarks                                                  	& type                        		\\
\hline
\hline
\multicolumn{8}{|c|}{direct integral-penalty finite element methods}                                                                                                                                                                                                                                    \\
\hline
\cite{Balakrishnan68}  	& not analyzed 		& 														& explicit ODE, $y \in L^2$ 			& $L^\infty$       			&    						& no path constraints                 						& not a numerical scheme             \\
\hline
\cite{DeJulio70}  		& $F$ strictly conv. bounded below  & & no state constraints 				& uniquely determines $y$       	&    						& no inequality constraints                 				&                                	\\
\hline
\cite{Jones70}   		& $C^1$ &                          & explicit ODE                		&                   				&  							& no mixed path-constraints                       			&                   				\\
\hline
\cite{Hager90}   		& $C^1$ 			& $C^0$                                 				& explicit ODE                		&                   				& SOSC						& no path-constraints                       				&                   				\\
\hline
\cite{McKinney72}   	& $f,c$ L-cont in $y$	&      		& explicit ODE, $y\in L^\infty \cap C^0$	& $L^\infty$           		&    						& no path-constraints                       				&                   				\\
\hline
\hline
\multicolumn{8}{|c|}{direct finite element collocation}                                                                                                                                                                                                                                               	\\
\hline
\cite{Maurer}   		& $C^1$   			& L-cont.                                          		& $C^1$                           	& $C^1$                           	& LICQ, SOSC 				& hom. cond., surj. lin. eq. constr.                       	& explicit Euler collocation               	\\
\hline
\cite{MR4046772}   		& $C^1$ 			& $C^1$                                            		& $W^{2,\infty}$                	& $W^{1,\infty}$                	&  							& hom. cond., surj. lin. eq. constr.                       	& implicit Euler collocation                  \\
\hline
\cite{kameswaran_biegler_2008}   			& suff. smooth & suff. smooth                         	& suff. smooth                  	& suff. smooth                    	& LICQ, SOSC				&                                                          	& hp collocation                          \\
\hline
\cite{RaoHager}   		& L-cont.      		& L-cont.                                          		& $H^2$                           	& $H^2$                           	& SOSC       				& $x^\star$ locally unique               					& hp coll.                          \\
\hline
\hline
\multicolumn{8}{|c|}{direct pseudo-spectral collocation}                                                                                                                                                                                                                                                \\
\hline
\cite{gongEtAl:2006}   	&              		&                                                   	& $C^{1}$ 							& $C^{1}$                           &            				& $\dot{x}_h$ converges uniformly      						& for feedback-linearizable systems \\
\hline
\cite{kangEtAl:2007}   &              		&                                                    	& $W^{1,\infty}$ 					& $W^{1,\infty}$ 					&            				&                                                          	& extends \cite{gongEtAl:2006}    	\\
\hline
\cite{gongEtAl:2008}   & $C^1$   & L-cont.                                                                      & $C^1$                             & $C^0$                           &            &                                                          & extends \cite{gongEtAl:2006},\cite{kangEtAl:2007}                    \\
\hline
\hline
\multicolumn{8}{|c|}{Runge-Kutta methods}                                                                                                                                                                                                                                               \\
\hline
\cite{SchwartzPolak:1996}   & L-cont.              &                                                                             & $W^{1,\infty}$       & $W^{1,\infty}$                                     &            & explicit ODE & mesh-uniform RKM        \\
\hline
\cite{Hager2000}   & $C^2$   &                                                                              & $W^{2,\infty}$ & $C^1$                           & SOSC       & no path constraints                                                         &     \\
\hline
\end{tabular}
\end{adjustbox}
\end{table}

	A convergence proof for the explicit Euler method is given in \cite{Maurer}, relying on local differentiability and Lipschitz properties of the problem-defining functions, smoothness of the optimal control solution, a homogeneity condition on the active constraints, surjectivity of linearized equality constraints, and a  coercivity assumption. These conditions are sophisticated, difficult to understand, and hard to verify and ensure by construction.
	The conditions assert the linear independence constraint qualification (LICQ) and the second-order sufficient condition (SOSC) of the NLP from the discretization, implying a unique and convergent local solution for the optimality conditions of the collocation method. This is why convergence proofs for other collocation schemes make similar assumptions.
	
	Implicit Euler for DAE constraints of index 2 is analyzed in \cite{MR4046772}. The authors prove convergence of critical points of NLP to critical points of \refOCP subject to assumptions on smoothness and local uniqueness.
	
	An $\ell_1$-penalty high-order collocation method is presented in \cite{Kameswaran} that appears stable for the problem of\cite{AlyChan} in Mayer form with path constraints. The result is experimental. Convergence results for a similar method are given in \cite{kameswaran_biegler_2008} for problems without path constraints. For higher order of convergence, problem-defining functions as well as the optimal control solution (including co-states) must satisfy additional smoothness properties.
	
	A high-order convergence result is presented in~\cite{RaoHager} when assuming sufficiently smooth, locally unique optimal control solutions. In contrast to \cite{kameswaran_biegler_2008}, they only require Lipschitz continuity of the problem-defining functions and their first derivatives.
	
	A convergence result for a pseudospectral method for the control of constrained feedback linearizable systems is given in~\cite{gongEtAl:2006}, relying on the strong assumption that the derivatives of the interpolating polynomial for the state converges uniformly to a continuous function. The assumption was relaxed in~\cite{kangEtAl:2007} to allow for discontinuous optimal inputs under piecewise differentiability assumptions.
	
	The assumption on feedback linearizable systems in~\cite{gongEtAl:2006,kangEtAl:2007} was relaxed in~\cite{gongEtAl:2008} to allow for more general nonlinear ordinary differential equations in trade for stronger requirements on the smoothness of the problem-defining functions and optimal control solution.
	
	Convergence proofs for shooting methods for quadratic regulator problems with linear dynamics and linear inequality constraints are presented in~\cite{pannochiaEtAl:2010,pannochiaEtAl:2015,yuzEtAl:2005}, based on matrix exponentials.
	
	High-order convergence proofs for explicit fixed-step size Runge-Kutta methods are given in~\cite{SchwartzPolak:1996}, extending results from~\cite{PolakBook1997} for the Euler methods. The analysis is limited to explicit time-invariant ODE constraints and norm-constrained inputs.
	
	Convergence results for Runge-Kutta methods for unconstrained optimal control problems are available in \cite{Hager2000}. In order to remove the input boundedness requirement, boundedness assumptions on the states and their derivatives are made, in conjunction with smoothness assumptions on the problem-defining functions and the optimal control solution.
	
	\subsubsection{Indirect Methods}
	Indirect methods have also been widely studied for computing the solution of~\refOCP~\cite{rao_asurvey,Betts2nd}. They use calculation of variations to determine the optimality conditions for the optimal arcs.
	Indirect methods have been less successful in practice than direct methods for these reasons:
	Firstly, the resulting optimality system, namely a Hamiltonian boundary-value problem (HBVP), is difficult to solve numerically. There are robustness issues when no accurate initial guess for the solution of the HBVP is given \cite{BrysonWildly}. Secondly, for complicated problems it is difficult or impossible to determine the optimality conditions to begin with, hence there are problems which cannot be solved by indirect methods \cite{Boehme2017}. Even when possible, state-constrained problems require estimates of the optimal switching structure. Workarounds for state constraints, such as by saturation functions, have been proposed in~\cite{WANG20174070}. Thirdly, for singular-arc problems, optimality conditions of higher order need to be used in order to determine a solution \cite[Sec.~1.4]{LAMNABHILAGARRIGUE1987173}. Finally, the optimality conditions can be ill-posed, for example when the co-state solution is non-unique. This always happens when path constraints are linearly dependent, because a dual solution is then non-unique.

	\subsection{Notation}
	\label{sec:assumptions}
	
	Let $-\infty<t_0 < t_E < \infty$ and the $M \in \N$ points
	$
	t_k \in \overline{\Omega}$,
	$\forall k\in \lbrace 1,2,\ldots,M \rbrace$. $\overline{\Omega}$ denotes the closure of $\Omega$.
	The functions
	$f : \R^{n_y} \times \R^{n_y} \times \R^{n_z} \times \Omega \rightarrow \R$,
	$c : \R^{n_y} \times \R^{n_y} \times \R^{n_z} \times \Omega \rightarrow \R^\nc$,
	$b: \R^{n_y} \times \R^{n_y} \times \ldots \times \R^{n_y} \rightarrow \R^\nb$.
	The function $y:\overline{\Omega} \rightarrow \R^{n_y}, t \mapsto y(t)$  and~$z:\overline{\Omega} \rightarrow \R^{n_z}, t \mapsto z(t)$.
	Given an interval $\Omega\subset \R$,  let $|\Omega|:=\int_\Omega 1\,\mathrm{d}t$.
	We use Big-$\cO$ notation to analyze a function's behaviour close to zero, i.e.\ function  $\phi(\xi)=\cO(\gamma(\xi))$ if and only if $\exists C>0$ and $\xi_0>0$ such that $\phi(\xi)\leq C \gamma(\xi)$ when $0<\xi<\xi_0$. The vector $\be:=[1\,\cdots\,1]\t$ with appropriate size.
	
	For notational convenience, we define the function
	\[
	x:=(y,z):\overline{\Omega} \rightarrow \R^{\nx},
	\]
	where $\nx:=n_y+n_z$.
	The solution space of $x$ is the Hilbert space
	\begin{align*}
		\cX := \left(H^1\left(\Omega\right)\right)^{n_y} \times \left(L^2\left(\Omega\right)\right)^{n_z}
	\end{align*}
	with scalar product
	\begin{align}
		    \langle(y,z),(v,w)\rangle_\cX := \sum_{j=1}^{n_y} \langle y_{[j]},v_{[j]} \rangle_{H^1(\Omega)} + \sum_{j=1}^{n_z} \langle z_{[j]},w_{[j]} \rangle_{L^2(\Omega)} \label{eqn:ScalarProd}
	\end{align}
	and induced norm $\|x\|_\cX := \sqrt{\langle x,x \rangle_\cX}$, where $\phi_{[j]}$ denotes the $j^\text{th}$ component of a function $\phi$. The Sobolev space $H^1(\Omega):=W^{1,2}(\Omega)$ and Lebesgue space $L^2(\Omega)$ with their respective scalar products are defined as in~\cite[Thm~3.6]{Adams}.
	The weak derivative of $y$ is denoted by
	$\dot{y} := {\mathrm{d}y}/{\mathrm{d}t}$.
	
	Recall the embedding $H^1(\Omega) \hookrightarrow \cC^0(\overline{\Omega})$, where $\cC^0(\overline{\Omega})$ denotes the space of continuous functions over $\overline{\Omega}$~\cite[Thm~5.4, part II, eqn~10]{Adams}. Hence, by requiring that $y\in \left(H^1\left(\Omega\right)\right)^{n_y}$ it follows that $y$ is continuous.
	In contrast, though $\dot{y}$ and $z$ are in $L^2(\Omega)$, they may be discontinuous.
	
	\subsection{Assumptions}
	In order to prove convergence, we make the following assumptions on \refOCP:
	\begin{enumerate}[\text{(A.}1\text{)}]
		\item 
		\refOCP has at least one global minimizer $x^\star$.
		\item 
		$\|c(\dot{y}(t),y(t),z(t),t)\|_1,\,$ and
		$\|b(y(t_1),\dots,y(t_M))\|_1$ are bounded for all arguments $x \in \cX$ within $z\geq 0$, $t \in \Omega$. $F(x)$ is bounded below for all arguments $x \in \cX$ within $z\geq 0$.
		\item 
		$f,\ c,\ b$ are globally Lipschitz continuous in all arguments except $t$.
		\item 
		The two solutions $x^\star_\omega,x^\star_{\omega,\tau}$ related to~$x^\star$, defined in Section~\ref{sec:reform}, are  bounded in terms of $\|z\|_{L^\infty(\Omega)}$ and $\|x\|_\cX$. Also, $\|x^\star\|_\cX$ is bounded.
		\item The related solution $x^\star_{\omega,\tau}$ can be approximated to an order of at least 1/2 using piecewise polynomials; formalized in \eqref{eqn:InfBound} below.
	\end{enumerate}
	Similar assumptions are  implicit or explicit in most of the literature. A discussion of these assumptions is appropriate:
	\begin{enumerate}[\text{(A.}1\text{)}]
	\item is just to avoid infeasible problems.
	
	\item The assumption on $b,c$ can be enforced by construction via lower and upper limits w.l.o.g.\ because they are (approximately) zero at the (numerical) solution. Boundedness below for $F$ is arguably mild when/since $\|x^\star\|_\cX,\|x_h\|_\cX$ are bounded: For minimum-time problems and positive semi-definite objectives this holds naturally. In many contexts, a lower bound can be given. The assumptions on $b,c$ have been made just to simplify the proof of a Lipschitz property and because they mean no practical restriction anyways. The boundedness assumption on $F$ is made to avoid unbounded problems.
	
	\item can be enforced. Functions that are not Lipschitz continuous, e.g.\ the square-root or Heaviside function, can be made so by replacing them with smoothed  functions, e.g.\ via a suitable mollifier. Smoothing is a common practice to ensure the derivatives used in a nonlinear optimization algorithm (e.g.\ IPOPT \cite{IPOPT}) are globally well-defined. The assumption has been made to prove a Lipschitz property of a penalty-barrier functional. Actually this property is only needed in a local neighborhood of the numerical optimal control solution, but for ease of notation we opted for global assumptions.
	
	\item can be ensured as shown in Remark~\ref{sec:ForceBoundedness} in Section~\ref{sec:PenaltyBarrierProblem}. This assumption effectively rules out the possibility of solutions with finite escape time. The assumption has been incorporated because restriction of a solution into a box means little practical restriction but significantly shortens convergence proofs due to boundedness.
	
	\item is rather mild, as discussed in Section~\ref{sec:InterpolationError} and illustrated in Appendix~\ref{app:3}. All finite-element methods based on piecewise polynomials make similar assumptions, implicitly or explicitly. The assumption is only used for the rate-of-convergence analysis. The assumption is unavoidable, since otherwise a solution $x^\star$ could exist that cannot be approximated to a certain order.
	\end{enumerate}
	The assumptions are not necessary but sufficient. Suppose that we have found a numerical solution. It is not of relevance to the numerical method whether the assumptions	hold outside of an open neighborhood of this solution. However, the proofs below would become considerably more lengthy with  local assumptions.  We outline in Section~\ref{sec:local} how our global analysis can be used to show local convergence under local assumptions. Hence, for the same reasons as in some of the literature, we opted for global assumptions. Our analysis is not restrictive in the sense that it imposes global requirements. 
	
	We do not make any further assumptions anywhere in the paper. Hence, our assumptions are milder than the methods discussed  in Section~\ref{sec:Lit2}, which can be technical and difficult to verify or enforce.
	
	\subsection{Outline}
	Section~\ref{sec:reform} introduces a reformulation of \refOCP as an unconstrained problem. Section~\ref{sec:FEM} presents the Finite Element Method in order to formulate a finite-dimensional unconstrained optimization problem. The main result of the paper is Theorem~\ref{thm:order}, which shows that solutions of the finite-dimensional optimization problem converge to solutions of~\refOCP with a guarantee on the order of convergence. Section~\ref{sec:solvingNLP} discusses how one could compute a solution using NLP solvers.
    Section~\ref{sec:numerical} presents numerical results which validate that our method converges for 
	difficult problems, whereas certain collocation methods can fail in some cases. Conclusions are drawn in Section~\ref{sec:conclusions}.

	\section{Reformulation as an Unconstrained Problem}\label{sec:reform}
	{\emergencystretch1em
	The reformulation of~\refOCP into an unconstrained problem is achieved in two steps. First, we introduce penalties for the 
	equality constraints. We then add logarithmic barriers for the inequality constraints.
	The resulting penalty-barrier functional will be treated numerically in Section~\ref{sec:FEM}.\par}
	
	Before proceeding, we note that boundedness and Lipschitz-continuity of $F$ and $r$ in \eqref{eqn:def:F}--\eqref{eqn:def:r} follow from (A.2)--(A.3).
	\begin{lem}[Boundedness and Lipschitz-continuity of $F$ and $r$]\label{lem:BoundLipschitz_Fr}
		$F$ is bounded below. $r$ is bounded. $F,r$ are Lipschitz continuous in $x$ with respect to
		$\|\cdot\|_\cX$. Furthermore, $F,r$ are Lipschitz continuous in $z$ with respect to the norm $\|\cdot\|_{L^1(\Omega)}$.
	\end{lem}
	The proof is given in Appendix~\ref{sec:Appendix_ProofLemma1}.
	
	We  bound the Lipschitz constants (i.e., with respect to both $\|x\|_\cX$ and $\|z\|_{L^1(\Omega)}$) with $L_F\geq 2$ for $F$ and with $L_r\geq 2$ for $r$.
	
	\subsection{Penalty Form}
	We introduce the \textit{penalty problem}
	\begin{equation}
	\tag{PP}
	\label{eqn:PP}
	\text{Find } x_\omega^\star \in \operatornamewithlimits{arg\, min}_{x \in \cX}
	F_\omega(x) \text{ s.t.\ } z(t)\geq 0 \text{ f.a.e.\ } t \in \Omega ,
	\end{equation}
	where $F_\omega(x) := F(x) + \frac{1}{2 \cdot \omega} \cdot r(x)$ and a small penalty parameter $\omega \in (0,1)$. Note that $F_\omega$ is Lipschitz continuous with constant
	\begin{align}
	    L_\omega := \max\left\lbrace L_F + \frac{L_r}{2 \omega}\,,\,L_f+\frac{L_c}{2 \omega}\, \|c\|_1\right\rbrace\,, \label{eqn:LipLw}
   \end{align}
	with $L_f,L_c$ the Lipschitz-constants of $f,c$ and $\|c\|_1$ is the upper bound on the 1-norm of $c$, as asserted by (A.2). We show that $\varepsilon$-optimal solutions of \eqref{eqn:PP} solve~\refOCP in a tolerance-accurate way.
	
	\begin{prop}[Penalty Solution]\label{thm:PenaltySolution}
	Let $\varepsilon\geq 0$. Consider an $\varepsilon$-optimal solution $x^\varepsilon_\omega$ to \eqref{eqn:PP}, i.e.
		$$ F_\omega(x^\varepsilon_\omega) \leq F_\omega(x^\star_\omega) + \varepsilon\text{ and } z^\varepsilon_\omega(t)\geq 0\ \text{f.a.e.}\ t\in\Omega\,. $$
		If we define $C_r := F(x^\star) - \operatornamewithlimits{ess\,min}_{x\in\cX, z\geq 0}F(x)$, then
$F(x^\varepsilon_\omega) \leq F(x^\star) + \varepsilon$, $r(x^\varepsilon_\omega) \leq \omega \cdot (C_r+\varepsilon)$.
	\end{prop}
	\begin{proof}
		$x^\star,\,x^\star_\omega,\,x^\varepsilon_\omega$ are all feasible for \eqref{eqn:PP}, but $x^\star_\omega$ is optimal and $x^\varepsilon_\omega$ is $\varepsilon$-optimal. Thus,
		\begin{align}
		F(x^\varepsilon_\omega) \leq F_\omega(x^\varepsilon_\omega) \leq 
	F_\omega(x^\star) + \varepsilon \leq F(x^\star) + \varepsilon
			\label{eqn:aux:prop1}
		\end{align}
		From this follows $F(x^\varepsilon_\omega)\leq F(x^\star)+\varepsilon$ because $r(x^\varepsilon_\omega)\geq 0$ and $r(x^\star)=0$ by (A.1). To show the second proposition, subtract $F(x^\varepsilon_\omega)$ from \eqref{eqn:aux:prop1}. Then it follows that  $1/(2\cdot\omega) \cdot r(x^\varepsilon_\omega) \leq F(x^\star) - F(x^\varepsilon_\omega) + \varepsilon \leq C_r + \varepsilon$\,. Multiplication of this inequality with $2\cdot\omega$ shows the result. Boundedness of $C_r$ follows from Lemma~\ref{lem:BoundLipschitz_Fr}.
	\end{proof}
	
	This result implies that for an $\varepsilon$-optimal solution  to~\eqref{eqn:PP} the optimality gap to~\refOCP is less than $\varepsilon$ and that the feasibility residual can be made arbitrarily small by choosing the  parameter~$\omega$ to be sufficiently small.

	\subsection{Penalty-Barrier Form}\label{sec:PenaltyBarrierProblem}
	We reformulate \eqref{eqn:PP} once more in order to remove the inequality constraints. We do so using logarithmic barriers. Consider the \textit{penalty-barrier problem}
	\begin{equation}
	\tag{PBP}
	\label{eqn:PBP}
	\text{Find }x^\star_{\omega,\tau} \in \operatornamewithlimits{arg\, min}_{x \in \cX}
	F_{\omega,\tau}(x) := F_\omega(x) + \tau \cdot \Gamma(x),
	\end{equation}
	where the barrier parameter $\tau \in (0,\omega]$ and $\Gamma$ is defined in~\eqref{eqn:def:Gamma}.
	
	We have introduced $\Gamma$ in order to keep $z^\star_{\omega,\tau}$ feasible with respect to \eqref{eqn:positive}. Recall that $L^2(\Omega)$ contains functions that have poles. So the following result is to ensure that $\Gamma$ actually fulfills its purpose.
	\begin{lem}[Strict Interiorness]\label{lem:StrictInteriorness}
		$$ 	z^\star_{\omega,\tau}(t) \geq \frac{\tau}{L_\omega}\cdot \qquad \text{f.a.e.\ } t \in	\Omega.$$
	\end{lem}
	\begin{proof}
	    At the minimizer $x^\star_{\omega,\tau}$, the functional $F_{\omega,\tau}$ can be expressed in a single component $z_{[j]}$ as
	    $$  \int_\Omega \Big(\,q\big( z_{[j]}(t),t \big) - \tau \cdot \log\big(\,z_{[j]}(t)\,\big)\,\Big)\,\mathrm{d}t\,,    $$
	    where $q$ is Lipschitz-continuous with a constant $L_q \leq L_\omega$ (cf. right argument in the max-expression \eqref{eqn:LipLw} and compare to \eqref{eqn:LipQuadPenaltyFunc} in the proof of Lemma~\ref{lem:BoundLipschitz_Fr}).
	    From the Euler-Lagrange equation it follows for $z^\star_{\omega,\tau}$ that
	    $$  \frac{\partial q}{\partial z_{[j]}}q\big( z_{[j]}(t) , t \big) - \frac{\tau}{z_{[j]}(t)} = 0\qquad \text{f.a.e.\ }t \in \Omega\,.    $$
	    The value of $z_{[j]}(t)$ gets closer to zero when the first term grows. However, that term is bounded by the Lipschitz constant. Hence, in the worst case
	    \[  z_{[j]}(t) \geq \frac{\tau}{L_q} \geq \frac{\tau}{L_\omega}\qquad \text{f.a.e.\ }t \in \Omega\,.\]
	\end{proof}
	
	We will need the following  operators:
	\begin{defn}[Interior Push]
		Given $x \in \cX$,  define $\bar{x}$ and $\check{x}$ as a modified  $x$ whose components $z$ have been pushed by an amount into the interior if they are close to zero:
		\begin{align*}
			\bar{z}_{[j]}(t):=\max\left\lbrace z_{[j]}(t),{\tau}/{L_\omega}\right\rbrace\,,\qquad \check{z}_{[j]}(t):=\max\left\lbrace z_{[j]}(t),{\tau}/(2\cdot L_\omega)\right\rbrace
		\end{align*}
		for all $j \in \lbrace 1,2,\ldots,n_z\rbrace$ and $t \in \overline{\Omega}$.
	\end{defn}
	\noindent
	Note that $\bar{x} \in \cX$ and that $x^\star_{\omega,\tau}=\bar{x}^\star_{\omega,\tau}$ from Lemma~\ref{lem:StrictInteriorness}.
	
	Using the interior push, we show below that $x^\star_{\omega,\tau}$ is $\varepsilon$-optimal for~\eqref{eqn:PP}. Our result uses a small arbitrary fixed number $0<\zeta\ll 1$.
	\begin{prop}[Penalty-Barrier Solution]\label{prop:PenaltyBarrierSolution}
		If (A.4) holds, then
		$$|F_\omega(x^\star_{\omega,\tau}) - F_\omega(x^\star_\omega)| = \cO\left(\tau^{1-\zeta}\right).$$
	\end{prop}
	\begin{proof}
		From the definition of the bar operator, we can use the bound
		\begin{align*}
			&\|x^\star_\omega - \bar{x}^\star_{\omega}\|_\cX = \|z^\star_{\omega}-\bar{z}^\star_{\omega}\|_{L^2(\Omega)}
			=\sqrt{\int_\Omega \|z^\star_{\omega}-\overline{z}^\star_{\omega}\|_2^2\,\mathrm{d}t}\\
			\leq& \operatornamewithlimits{max}_j\sqrt{|\Omega| \cdot n_z \cdot \|z^\star_{\omega\,[j]}-\overline{z}^\star_{\omega\,[j]}\|^2_{L^\infty(\Omega)}}\leq n_z \cdot \sqrt{|\Omega|} \cdot \frac{\tau}{L_\omega},
		\end{align*}
		together with the facts that $x^\star_{\omega,\tau}= \bar{x}^\star_{\omega,\tau}$ and $F_\omega$ is Lipschitz continuous, to get
		\begin{align*}
			0 &\leq
			F_\omega(x^\star_{\omega,\tau}) - F_\omega(x^\star_{\omega}) \leq  F_\omega(\bar{x}^\star_{\omega,\tau}) - F_\omega(\bar{x}_\omega^\star) + L_\omega \cdot \|x^\star_{\omega}-\bar{x}^\star_{\omega}\|_\cX\\
			& \leq \underbrace{F_\omega(\bar{x}^\star_{\omega,\tau}) - F_{\omega,\tau}(\bar{x}^\star_{\omega,\tau})}_{=-\tau\cdot\Gamma(\bar{x}^\star_{\omega,\tau})} + F_{\omega,\tau}(\bar{x}^\star_{\omega,\tau})\\
			&\phantom{ \leq }-\Big(\underbrace{F_\omega(\bar{x}^\star_{\omega}) - F_{\omega,\tau}(\bar{x}^\star_{\omega})}_{=-\tau \cdot \Gamma(\bar{x}^\star_{\omega})} + F_{\omega,\tau}(\bar{x}^\star_{\omega})\Big)+ L_\omega \cdot n_z \cdot \sqrt{|\Omega|} \cdot \frac{\tau}{L_\omega}\\
			& \leq F_{\omega,\tau}(\bar{x}^\star_{\omega,\tau})-F_{\omega,\tau}(\bar{x}^\star_{\omega}) +|\tau \cdot \Gamma(\bar{x}^\star_{\omega,\tau})| + |\tau \cdot \Gamma(\bar{x}^\star_\omega)| + n_z \cdot \sqrt{|\Omega|} \cdot \tau .
		\end{align*}
		We use $|\tau \cdot \Gamma(\bar{x}^\star_{\omega,\tau})|+|\tau \cdot \Gamma(\bar{x}^\star_\omega)|=\cO(\tau^{1-\zeta})$, as per Lemma~\ref{lem:Gamma} in Appendix~\ref{sec:Appendix_BarrierFunctioProperties}, to obtain the result from
		\begin{align*}
			\hspace{1mm}F_\omega(x^\star_{\omega,\tau}) - F_\omega(x^\star_{\omega})
			\leq \underbrace{F_{\omega,\tau}(\bar{x}^\star_{\omega,\tau})-F_{\omega,\tau}(\bar{x}^\star_{\omega})}_{\leq 0} +	\cO\left(\tau^{1-\zeta}\right)+ n_z \cdot \sqrt{|\Omega|} \cdot \tau
		\end{align*}
		The under-braced term is bounded above by zero because $\bar{x}^\star_{\omega,\tau}= x^\star_{\omega,\tau}$ is a minimizer of $F_{\omega,\tau}$.
	\end{proof}
	
	\begin{rem}	\label{sec:ForceBoundedness}
		Lemma~\ref{lem:Gamma} in the proof of Prop.~\ref{prop:PenaltyBarrierSolution} needs (A.4), i.e.
		\[\|z^\star_{\omega}\|_{L^\infty(\Omega)},\|z^\star_{\omega,\tau}\|_{L^\infty(\Omega)} = \cO(1).\]
		Note that the assumption can be enforced. For example, the path constraints
		$$ 	z_{[1]}(t)\geq 0,\quad z_{[2]}(t)\geq 0,\quad z_{[1]}(t)+z_{[2]}(t)=const$$
		lead to
		$  \|z_{[j]}\|_{L^\infty(\Omega)} \leq const$ \mbox{for $j=1,2\,.$}
		Constraints like these arise when variables have simple upper and lower  bounds before being transformed 
		into \refOCP.
		
		Similarly, boundedness of $\|x\|_\cX$ can be enforced. To this end, introduce box constraints for each component of $\dot{y},y,z$, before transcribing into the form \refOCP.
	\end{rem}
	
	\section{Finite Element Method}\label{sec:FEM}
	Our method constructs an approximate finite element solution $x^\epsilon_h$ by solving the unconstrained problem \eqref{eqn:PBP} computationally in a finite element space \mbox{$\cX_{h,p} \subset \cX$,} using an NLP solver.
	
	We introduce a suitable finite element space and show a stability result in this space. Eventually, we prove convergence of the finite element solution to  solutions of \eqref{eqn:PBP} and \refOCP. 
	
	\subsection{Definition of the Finite Element Space}
	Let the mesh parameter $h \in (0,|\Omega|]$.
	The set	$\cT_h$ is called a \emph{mesh} and consists of open intervals $T \subset \Omega$ that satisfy the usual conditions~\cite[Chap.~2]{Ciarlet1978}:
	\begin{enumerate}
	    \item Disjunction: $T_1\cap T_2 = \emptyset$, for all  distinct  $T_1,T_2 \in \cT_h$.
	    \item Coverage: $\bigcup_{T \in \cT_h} \overline{T} = \overline{\Omega}$.
	    \item Resolution: $\max_{T \in \cT_h} |T| = h$.
	    \item Quasi-uniformity: $\min_{T_1,T_2 \in \cT_h}\frac{|T_1|}{|T_2|} \geq \vartheta > 0$. The constant $\vartheta$ must not depend on $h$ and $1/\vartheta = \cO(1)$.
	\end{enumerate}
	
	We write $\cP_p({T})$ for the space of functions that are polynomials of degree $\leq p \in \N_0$ on interval ${T}$. Our finite element space is then given as
	\begin{align*}
		\cX_{h,p} := \left\{x:\overline{\Omega}\rightarrow\R^{\nx}  \mid y \in \cC^0(\overline{\Omega}), x \in \cP_p(T)^{n_x} \ \forall T \in \cT_h \right\}.
	\end{align*}
	$\cX_{h,p}\subset \cX$ is a Hilbert space with scalar product $\langle\cdot,\cdot\rangle_\cX$. 
	
\begin{figure}[tb]
		\centering
		\includegraphics[width=0.6\columnwidth]{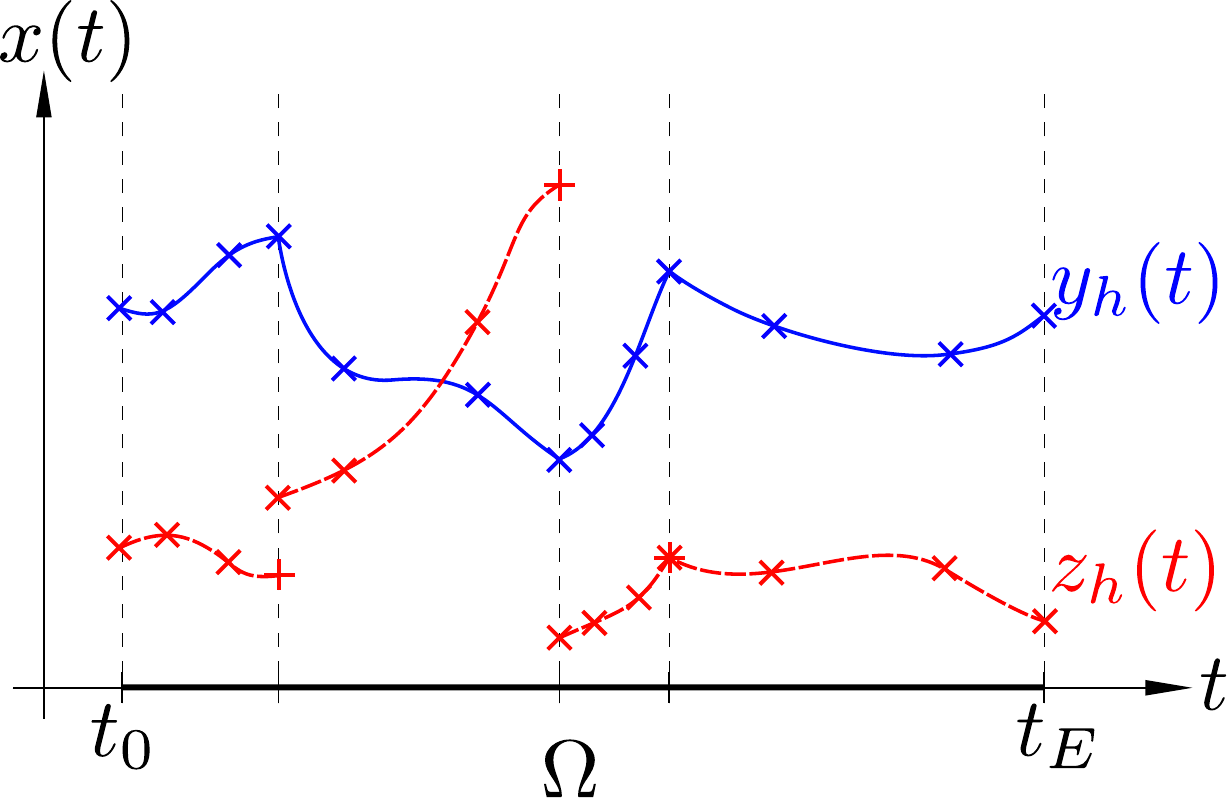}
		\caption{Continuous and discontinuous piecewise polynomial finite element functions $y_h,z_h$ on a mesh $\cT_h$ of four intervals.}
		\label{fig:finiteelementsyz}
	\end{figure}

	Note that if $(y,z)\in\cX_{h,p}$, then $y$ is continuous but~$\dot{y}$ and~$z$ can be discontinuous. Figure~\ref{fig:finiteelementsyz} illustrates two functions $(y_h,z_h) \in \cX_{h,p}$ with $\times$ and $+$ for their nodal basis, to identify them with a finite-dimensional vector. 
	
	\subsection{Discrete Penalty-Barrier Problem}
	We state the \textit{discrete penalty-barrier problem}
	as
	\begin{equation*}
		\text{Find }x^\star_h \in \operatornamewithlimits{arg\,min}_{x \in \cX^{\omega,\tau}_{h,p}}\, F_{\omega,\tau}(x)
		\tag{PBP\textsubscript{h}}
        \label{eqn:PBP_h}
	\end{equation*}
	with the space $\cX^{\omega,\tau}_{h,p} := \left\lbrace\,x \in \cX_{h,p} \ \Big\vert \ z(t)\geq \frac{\tau}{2 \cdot L_\omega}\cdot\be \text{ f.a.e.\ }t \in\Omega\right\rbrace$.
	
	Note that Lemma~\ref{lem:StrictInteriorness} is valid only for solutions to~\eqref{eqn:PBP}, whereas below we will consider sub-optimal solutions to~\eqref{eqn:PBP_h}. Hence, we cannot guarantee that these sub-optimal solutions will satisfy $z(t)\geq {\tau}/(1 \cdot L_\omega)\cdot\be$. The looser constraint $z(t)\geq {\tau}/(2 \cdot L_\omega)\cdot\be$ in the definition above will be used in the proof of Theorem~\ref{thm:ConvOCP}.
	
	In a practical implementation,
	we neglect these additional constraints. This is reasonable when solving the NLP with interior-point methods, since they keep the numerical solution strictly interior with a distance to zero in the order of $\tau \gg \frac{\tau}{2 \cdot L_\omega}$.
	
	\subsection{Stability}
	The following result shows that two particular Lebesgue norms are equivalent in the above finite element space.
	\begin{lem}[Norm equivalence]\label{lem:NormEquivalence}
		If $x \in \cX_{h,p}$, then
		\begin{align*}
			\|x_{[j]}\|_{L^\infty(\Omega)} \leq \frac{p+1}{\sqrt{\vartheta \cdot h}} \cdot \|x\|_\cX\quad\forall j\in\lbrace 1,2,\ldots,\nx\rbrace.
		\end{align*}
	\end{lem}
	\begin{proof}
		We can bound
		$\|x_{[j]}\|_{L^\infty(\Omega)}
		\leq\operatornamewithlimits{max}_{T \in \cT_h} \|x_{[j]}\|_{L^\infty(T)}$. We now use \eqref{eqn:PropAppendix2} in Appendix~\ref{sec:Appendix_LebesgueIdentity}. Since $x_{[j]} \in \cP_p(T)$, it follows that
		\begin{align*}
			\operatornamewithlimits{max}_{T \in \cT_h}\|x_{[j]}\|_{L^\infty(T)}
			&\leq \operatornamewithlimits{max}_{T \in \cT_h} \frac{p+1}{\sqrt{|T|}} \cdot \|x_{[j]}\|_{L^2(T)}\leq\frac{p+1}{\sqrt{\vartheta \cdot h}} \cdot \|x_{[j]}\|_{L^2(\Omega)} \leq \frac{p+1}{\sqrt{\vartheta \cdot h}} \cdot \|x\|_\cX.
		\end{align*}
	\end{proof}
	
	Below, with the help of Lemma~\ref{lem:NormEquivalence}, we obtain a bound on the growth of $F_{\omega,\tau}$ in a neighborhood of a solution $x^\star_{\omega,\tau}$ to~\eqref{eqn:PBP} for elements in~$\cX_{h,p}$.
	\begin{prop}[Lipschitz continuity]\label{prop:Lcont}
		Let
		\begin{align*}
			\delta_{\omega,\tau,h} &:= \frac{\tau}{2 \cdot L_\omega} \cdot \frac{\sqrt{\vartheta \cdot h}}{p+1}\,,\qquad
			L_{\omega,\tau,h} := L_\omega + n_z \cdot |\Omega| \cdot 2 \cdot L_\omega \cdot \frac{p+1}{\sqrt{\vartheta \cdot h}}.
		\end{align*}
		Consider the spherical neighbourhood
		\begin{align*}
			\cB := \big\lbrace\,x \in \cX\ \big\vert\ \|x^\star_{\omega,\tau}-x\|_\cX \leq \delta_{\omega,\tau,h} \big\rbrace.
		\end{align*}
		The following holds $\forall x^\text{A},x^\text{B} \in \cB\cap\cX_{h,p}$:
		\begin{align*}
			|F_{\omega,\tau}(x^\text{A}) - F_{\omega,\tau}(x^\text{B})| \leq L_{\omega,\tau,h} \cdot \|x^\text{A} - x^\text{B}\|_\cX.
		\end{align*}
	\end{prop}
	\begin{proof}
		From Lemma~\ref{lem:StrictInteriorness} and Lemma~\ref{lem:NormEquivalence} follows:
		\begin{align*}
  			\operatornamewithlimits{ess\,inf}_{t \in \Omega} z_{[j]}(t)\geq&\underbrace{\operatornamewithlimits{ess\,inf}_{t \in \Omega} z^\star_{\omega,\tau,[j]}(t)}_{\geq \frac{\tau}{L_\omega}}
			- \underbrace{\|z^\star_{\omega,\tau,[j]}-z_{[j]}\|_{L^\infty(\Omega)}}_{\leq \frac{p+1}{\sqrt{\vartheta \cdot h}} \cdot \delta_{\omega,\tau,h}\leq \frac{\tau}{2 \cdot L_\omega}}\hspace{4mm} \forall x \in \cB \cap \cX_{h,p}
		\end{align*}
		Hence,
		\begin{align}
			\operatornamewithlimits{min}_{1 \leq j \leq n_z} \operatornamewithlimits{ess\,inf}_{t\in\Omega} z_{[j]}(t) \geq \frac{\tau}{2 \cdot L_\omega}\quad \forall x \in \cB \cap \cX_{h,p}.\label{eqn:aux1:Lcont}
		\end{align}
		From Lipschitz-continuity of $F_\omega$ we find
		\begin{align*}
			    &|F_{\omega,\tau}(x^\text{A})-F_{\omega,\tau}(x^\text{B})|\\
			\leq&|F_{\omega}(x^\text{A})-F_{\omega}(x^\text{B})| + \tau \cdot \sum_{j=1}^{n_z} \int_\Omega\,\left|\log\left(z^\text{A}_{[j]}(t)\right)-\log\left(z^\text{B}_{[j]}(t)\right)\right|\,\mathrm{d}t \\
			\leq&L_\omega \cdot \|x^\text{A}-x^\text{B}\|_\cX + \tau \cdot n_z \cdot |\Omega|\cdot\operatornamewithlimits{max}_{1\leq j \leq n_z}\ \operatornamewithlimits{ess\,sup}_{t \in \Omega} \left|\log\left(z^\text{A}_{[j]}(t)\right)-\log\left(z^\text{B}_{[j]}(t)\right)\right|.
		\end{align*}
		We know a lower bound for the arguments in the logarithms from~\eqref{eqn:aux1:Lcont}. Thus, the essential supremum term can be bounded with a Lipschitz result for the logarithm:
		\begin{align*}
			&\operatornamewithlimits{max}_{1\leq j \leq n_z}\ \operatornamewithlimits{ess\,sup}_{t \in \Omega} \left|\log\left(z^\text{A}_{[j]}(t)\right)-\log\left(z^\text{B}_{[j]}(t)\right)\right|\\
			&\quad \leq\operatornamewithlimits{max}_{1\leq j \leq n_z}\ \frac{1}{\, \frac{\tau}{2 \cdot L_\omega} \,} \cdot \|z_{[j]}^\text{A}-z_{[j]}^\text{B}\|_{L^\infty(\Omega)}
			\leq \frac{2 \cdot L_\omega}{\tau} \cdot\frac{p+1}{\sqrt{\vartheta \cdot h}} \cdot \|x^\text{A}-x^\text{B}\|_\cX,
		\end{align*}
		where the latter inequality is obtained using Lemma~\ref{lem:NormEquivalence}. \end{proof}
	
	\subsection{Interpolation Error}\label{sec:InterpolationError}
	In order to show high-order convergence results, it is imperative that the solution function can be represented with high accuracy in a finite element space. In the following we introduce a suitable assumption for this purpose.
	
	Motivated by the Bramble-Hilbert Lemma \cite{BrambleHilbert}, we make the assumption (A.5) that for a fixed chosen degree $p = \cO(1)$ there exists an $ \ell \in (0,\infty)$ such that
	\begin{align}
		\infexpr = \cO\big(h^{\ell+1/2}\big)\,.  \label{eqn:InfBound}
	\end{align}
	Notice that the best approximation $x_h$ is well-defined since $\cX_{h,p}$ is a Hilbert space with induced norm $\|\cdot\|_\cX$. In Appendix~\ref{app:3} we give two examples to demonstrate the mildness of assumption~\eqref{eqn:InfBound}.
	
	To clarify on the mildness of \eqref{eqn:InfBound}, consider the triangular inequality
	\begin{align}
	    \infexpr \leq \|x^\star_{\omega,\tau} - x^\star\|_\cX + \operatornamewithlimits{min}_{x_h \in \cX_{h,p}}\|x^\star-x_h\|_\cX\,,
	\end{align}
	where $x^\star$ is the global minimizer. Clearly, the second term converges under the approximability assumption of finite elements, hence could not be milder. The first term holds under several sufficient assumptions; for instance if $x^\star$ is unique because then convergence of feasibility residual and optimality gap will determine --at a convergence rate depending on the problem instance-- that unique solution. Due to round-off errors in computations on digital computers, for numerical methods the notion of well-posedness is imperative, hence must always be assumed, relating to how fast the first term converges as optimality gap and feasibility residual converge.
	
	For the remainder, we define $\nu:= \ell/2$, $\eta:=(1-\zeta)\cdot\nu$ with respect to $\ell,\zeta$. We choose $\tau = \cO(h^{\nu})$ and $\omega = \cO(h^{\eta}) $ with
	$h>0$ suitably small such that $0 < \tau \leq \omega < 1$.
	
	Following the assumption \eqref{eqn:InfBound}, the result below shows that the best approximation in the finite element space satisfies an approximation property.
	\begin{lem}[Finite Element Approximation Property]
		\label{lem:bestApproximation}
		If~\eqref{eqn:InfBound} holds and  $h>0$ is chosen sufficiently small, then
		\begin{equation}
		\infexpr \leq \delta_{\omega,\tau,h}. \label{eqn:SphereBound}
		\end{equation}
	\end{lem}
	\begin{proof}
		For $h>0$ sufficiently small it follows from $\ell> \nu + \eta$, that $h^{\ell+1/2} < h^{\nu+\eta+1/2}$. Hence,
		$$	\infexpr\leq const \cdot h^{\nu+\eta+1/2}	$$
		for some constant $const$. The result follows by noting that
		\begin{align*}	
			\delta_{\omega,\tau,h}
			& 
			\geq \frac{\tau}{\,\frac{L_r}{\omega}\,} \cdot \frac{\sqrt{\vartheta \cdot h}}{p+1}
			= \frac{\sqrt{\vartheta}}{L_r\cdot(p+1)} \cdot \tau \cdot \omega \cdot \sqrt{h} \geq const \cdot h^{\nu+\eta+1/2}.
		\end{align*}
	\end{proof}
	In other words, Lemma~\ref{lem:bestApproximation} says for $h>0$ sufficiently small it follows that $\cB \cap \cX_{h,p}\neq \emptyset$. This is because the minimizing argument of \eqref{eqn:SphereBound} is an element of $\cB$.
	
	\subsection{Optimality}
	We show that an $\epsilon$-optimal solution for \eqref{eqn:PBP_h} is an $\varepsilon$-optimal solution for \eqref{eqn:PBP}, where $\varepsilon \geq \epsilon$.
	\begin{thm}[Optimality of Unconstrained FEM Minimizer]\label{thm:conv}
		Let $\cB$ as in Proposition~\ref{prop:Lcont}, and $x^\epsilon_h$ an $\epsilon$-optimal solution for \eqref{eqn:PBP_h}, i.e.	
		\begin{align*}
			F_{\omega,\tau}(x^\epsilon_h) \leq F_{\omega,\tau}(x^\star_h) + \epsilon.
		\end{align*}
		If $\cB \cap \cX_{h,p}\neq \emptyset$, then $x^\epsilon_h$ satisfies:
		\begin{align*}
			F_{\omega,\tau}(x^\epsilon_h)
			\leq F_{\omega,\tau}(x^\star_{\omega,\tau})+ \epsilon+ L_{\omega,\tau,h} \cdot \infexpr .
		\end{align*}
	\end{thm}
	\begin{proof}
		Consider the unique finite element best approximation from \eqref{eqn:SphereBound}
		$$ 	\tilde{x}_h := \operatornamewithlimits{arg\,min}_{x_h \in\cX_{h,p}}\|x^\star_{\omega,\tau}-x_h\|_\cX\,.  $$
		Since $\cB\cap \cX_{h,p} \neq \emptyset$, it follows
		$\tilde{x}_h \in \cB \cap \cX_{h,p}$. Hence,
		$$ 	\tilde{x}_h = \operatornamewithlimits{arg\,min}_{x_h \in\cB\cap\cX_{h,p}}\|x^\star_{\omega,\tau}-x_h\|_\cX .  $$
		From \eqref{eqn:aux1:Lcont} we find $\cB \cap \cX_{h,p} \subset \cX_{h,p}^{\omega,\tau}$. Thus, $\tilde{x}_h \in \cX^{\omega,\tau}_{h,p}$. Hence,
		$$ 	\tilde{x}_h = \operatornamewithlimits{arg\,min}_{x_h \in\cX^{\omega,\tau}_{h,p}}\|x^\star_{\omega,\tau}-x_h\|_\cX\,.  $$
		
		Proposition~\ref{prop:Lcont} can be used to obtain 
		$	F_{\omega,\tau}(\tilde{x}_h) \leq F_{\omega,\tau}(x^\star_{\omega,\tau}) + L_{\omega,\tau,h} \cdot \|x^\star_{\omega,\tau}-\tilde{x}_h\|_\cX$. 
		Since $x^\epsilon_h$ is a global $\epsilon$-optimal minimizer of $F_{\omega,\tau}$ in $\cX^{\omega,\tau}_{h,p}$ and also $\tilde{x}_h$ lives in $\cX^{\omega,\tau}_{h,p}$, the optimalities must relate as
			$F_{\omega,\tau}(x^\epsilon_h) \leq F_{\omega,\tau}(\tilde{x}_h)+\epsilon$. 
		The result follows. 
	\end{proof}
	
	\subsection{Convergence}
	We obtain a bound for the optimality gap and feasibility residual of $x^\epsilon_h$.
	
	\begin{thm}[Convergence to~\refOCP]\label{thm:ConvOCP}
		Let $x^\epsilon_h$ be an $\epsilon$-optimal numerical solution to~\eqref{eqn:PBP_h}. If (A.4) holds, then $x^\epsilon_h$ satisfies
		\begin{align*}
			g_\mathrm{opt} = \cO\left(\tau^{1-\zeta}+\varepsilon_{h,p} \right),\
			r_\mathrm{feas} = \cO\big(\omega \cdot \left(1+ \tau^{1-\zeta} + \varepsilon_{h,p}\right)\big),
		\end{align*}
		where
		$$ 	\varepsilon_{h,p} := L_{\omega,\tau,h} \cdot \infexpr + \epsilon\,.	$$
	\end{thm}
	\begin{proof}
		From Theorem~\ref{thm:conv} we know \mbox{$F_{\omega,\tau}(x^\epsilon_h) \leq F_{\omega,\tau}(x^\star_{\omega,\tau}) + \varepsilon_{h,p}$}. This is equivalent to
		\begin{align*}
			&F_\omega(x^\epsilon_h) + \tau \cdot \Gamma(x^\epsilon_h) \leq F_\omega(x^\star_{\omega,\tau}) + \tau \cdot \Gamma(x^\star_{\omega,\tau}) + \varepsilon_{h,p} \\
			\Rightarrow\ \,
			&F_\omega(x^\epsilon_h) \leq F_\omega(x^\star_{\omega,\tau})+ \underbrace{|\tau \cdot \Gamma(x^\epsilon_h)| + |\tau \cdot \Gamma(x^\star_{\omega,\tau})|}_{(*)} + \varepsilon_{h,p}.
		\end{align*}
		Since $x^\epsilon_h \in \cX^{\omega,\tau}_{h,p}$, it follows that $z^\epsilon_h \geq \frac{\tau}{2 \cdot L_\omega}\cdot \be$
		and thus $x^\epsilon_h = \check{x}^\epsilon_h$.
		From Lemma~\ref{lem:StrictInteriorness} we know $x^\star_{\omega,\tau} = \bar{x}^\star_{\omega,\tau}$. Thus, we can apply Lemma~\ref{lem:Gamma} to bound $(*)$ with $\cO(\tau^{1-\zeta})$. Hence,
			$F_\omega(x^\epsilon_h) \leq F_\omega(x^\star_{\omega,\tau}) + \cO(\tau^{1-\zeta}) + \varepsilon_{h,p}$.
		Since, according to Proposition~\ref{prop:PenaltyBarrierSolution}, $x^\star_{\omega,\tau}$ is $\tilde{\varepsilon}$-optimal for~\eqref{eqn:PP}, where
		$\tilde{\varepsilon} = \cO(\tau^{1-\zeta}),$
		it follows that
		\begin{align*}
			F_\omega(x^\epsilon_h) \leq F_\omega(x^\star_{\omega}) + \underbrace{\cO(\tau^{1-\zeta}) + \varepsilon_{h,p}}_{=:\varepsilon}.
		\end{align*}
		In other words, $x^\epsilon_h$ is $\varepsilon$-optimal for \eqref{eqn:PP}. The result now follows from Proposition~\ref{thm:PenaltySolution}.
	\end{proof}
	
	
	Below, we translate the above theorem into an order-of-convergence result.
	\begin{thm}[Order of Convergence to~\refOCP]\label{thm:order}
		Consider $x^\epsilon_h$ with $\epsilon = \cO(h^{\ell-\eta})$. Then
			$g_\mathrm{opt} = \cO\left(h^\eta\right)$ and
			$r_\mathrm{feas} =\cO\left(h^\eta\right)$.
	\end{thm}
	\begin{proof}
		It holds that
		\begin{align*}
			L_{\omega,\tau,h}
			&= \,\,\,\,L_\omega\,\,\, + (n_z \cdot |\Omega| \cdot 2) \cdot\,\, L_\omega\,\,\,\, \cdot \left(\frac{p+1}{\sqrt{\vartheta \cdot h}}\right)\\
			&=\frac{\cO(1)}{\omega} +\,\,\,\,\,\,\,\,\,\, \cO(1) \,\,\,\,\,\,\,\,\,\cdot \frac{\cO(1)}{\omega} \cdot \,\,\,\,\frac{\cO(1)}{\sqrt{h}}\\
			&=\cO\big(h^{-\eta} + h^{-\eta} \cdot h^{-1/2} \big) = \cO(h^{-\eta -1/2}).
		\end{align*}		
		From \eqref{eqn:InfBound}, we find
		\begin{align*}
			\varepsilon_{h,p} 	
			&= L_{\omega,\tau,h} \cdot \infexpr + \epsilon\\
			&= \cO\big(h^{-\eta -1/2}\big) \cdot \cO\big(h^{\ell+1/2}\big) + \cO\big(h^{\ell-\eta}\big)= \cO\big(h^{\ell-\eta}\big).
		\end{align*}
		Combining this with Theorem~\ref{thm:ConvOCP}, we find $x^\epsilon_h$ satisfies
		\begin{align*}
			g_\text{opt}
			=& \cO\big(\tau^{1-\zeta} + \varepsilon_{h,p}\big)
			=\cO\big(h^{\nu \cdot (1-\zeta)}+h^{\ell-\eta}\big)
			=\cO\big(h^{\min\left\lbrace \eta ,\, \ell-\eta \right\rbrace}\big),\\
			r_\text{feas}
			=& \cO\big(\omega \cdot (1+\tau^{(1-\zeta)\cdot\nu} + \varepsilon_{h,p})\big)
			= \cO\big(h^\eta + h^{\eta+(1-\zeta)\cdot\nu} + h^{\eta+\ell-\eta} \big)\\
			=& \cO\big( h^{\eta} + h^{2 \cdot\eta} + h^{\ell} \big)= \cO\big(h^{\min\lbrace \eta ,\, \ell \rbrace}\big).
		\end{align*}
		Note that $\ell > \ell-\eta>\eta$.
	\end{proof}
	Recall that $\eta = \ell/2 \cdot (1-\zeta)$, where $0<\zeta\ll 1$. If $\ell\approx p$ and $\eta\approx \ell/2$ it follows that $h^\eta \approx \sqrt{h^p}$.
	
	\subsection{Numerical Quadrature}
	When computing $x^\epsilon_h$, usually the integrals in $F$ and~$r$ cannot be evaluated exactly. In this case, one uses numerical quadrature and replaces $F_{\omega,\tau}$ with
	$F_{\omega,\tau,h} := F_h + \frac{1}{2 \cdot \omega} \cdot r_h + \tau \cdot \Gamma$. Since $\cX_{h,p}$ is a space of piecewise polynomials, $\Gamma$ can be integrated analytically. However, the analytic integral expressions become very complicated. This is why, for a practical method, one may also wish to use quadrature for $\Gamma$.
	
	If $F$ and $r$ have been replaced with quadrature approximations $F_h,\,r_h$, then it is sufficient that these approximations satisfy
	\begin{align}
		|F_{\omega,\tau,h}(x)-F_{\omega,\tau}(x)| \leq C_\text{quad} \cdot \frac{h^q}{\omega}\quad\quad \forall x \in \cX^{\omega,\tau}_{h,p},\label{eqn:QuadCond}
	\end{align}
	with bounded constant $C_\text{quad}$ and quadrature order $q\in\N$, to ensure that  the convergence theory holds. We discuss this further below.
    The constrraint \eqref{eqn:QuadCond} poses a consistency and stability condition.
	
	\paragraph{Consistency}
	There is a consistency condition in~\eqref{eqn:QuadCond} that relates to suitable values of $q$. In particular, if we want to ensure convergence of order $\cO(h^\eta)$, as presented in Theorem~\ref{thm:order}, then $q$ has to be sufficiently large.
	
	Consider the problem
	\begin{equation*}
		\tilde{x}^\star_{h} \in \operatornamewithlimits{arg\,min}_{x \in \cX_{h,p}^{\omega,\tau}}\ F_{\omega,\tau,h}(x).
	\end{equation*}
	Note that $\tilde{x}^\star_{h}$ is $\epsilon$-optimal for \eqref{eqn:PBP_h}, where from
	\begin{align*}
		F_{\omega,\tau}(\tilde{x}^\star_h)-C_\text{quad} \cdot \frac{h^q}{\omega} &\leq F_{\omega,\tau,h}(\tilde{x}^\star_h)
		\leq F_{\omega,\tau,h}(x^\star_h) \leq F_{\omega,\tau}(x^\star_h)+C_\text{quad} \cdot \frac{h^q}{\omega}
	\end{align*}
	it follows that $ 	\epsilon = \cO\left({h^q}/{\omega}\right) = \cO(h^{q-\eta}).  	$
	Hence, $\tilde{x}^\star_{h}$ satisfies the bounds for the optimality gap and feasibility residual presented in Theorem~\ref{thm:ConvOCP}.  We obtain the same order of convergence as in Theorem~\ref{thm:order} when maintaining $\epsilon = \cO(h^{\ell-\eta})$, i.e.\ choosing $q\geq \ell$.
	
	\paragraph{Stability}
	Beyond consistency, \eqref{eqn:QuadCond} poses a non-trivial stability condition. This is because the error bound must hold $\forall x \in \cX_{h,p}$. We show this with an example.
	
	Consider $\Omega=(0,1)$, $n_y=0$, $n_z=1$, and \mbox{$c(x):=\sin(\pi \cdot x)$.} The constraint forces $x(t)=0$. Clearly, $c$ and $\nabla c$ are bounded globally. Consider the uniform mesh $\cT_h:=\left\lbrace\,T_j \ \vert \ T_j = \left((j-1)\cdot h, j\cdot h\right), j=1,2,\ldots,1/h \right\rbrace$
	for $h \in 1/\N$, choose $p=1$ for the finite element degree, and Gauss-Legendre quadrature of order $q=3$, i.e.\ the mid-point rule quadrature scheme \cite{GaussQuad} of $n_q=1$ point per interval. Then, the finite element function $x_h$, defined as
	$x(t) := -{1}/{h} + {2}/{h} \cdot (t-j \cdot h)$ for $t \in T_j$ on each interval, yields the quadrature error
	\begin{align*} 	
		|r_h(x)-r(x)| = \Bigg|\underbrace{ h \cdot \sum_{j=1}^{1/h} \sin^2\big(\pi \cdot x(j\cdot h - h/2)\big)}_{= 0} - \underbrace{\int_{0}^{1}\sin^2\big(\pi \cdot x(t)\big)\mathrm{d}t}_{=0.5} \Bigg|,
	\end{align*}
	violating \eqref{eqn:QuadCond}. In contrast, using Gauss-Legendre quadrature of order $5$ (i.e.\ using $n_q=2$ quadrature points per interval) yields satisfaction of \eqref{eqn:QuadCond} with $q=5$.
	
	We see that in order to satisfy \eqref{eqn:QuadCond}, a suitable quadrature rule must take into account the polynomial degree $p$ of the finite element space and the nature of the nonlinearity of $c$. We clarify this using the notation $(\phi \circ \psi)(\cdot):=\phi(\psi(\cdot))$ for function compositions: If
	$$ 	\left(f + \frac{1}{2\cdot\omega}\cdot \|c\|_2^2 \right) \circ x \in \cP_d({T})^{n_x},\ \forall T\in\cT_h,\forall x \in \cX_{h,p}\cap\cB,$$
	for some $d \in \N$, i.e.\ the integrands of $F$ and~$r$ are polynomials in $t$, then $q \geq d$ is a sufficient order for exact quadrature. For a practical method, we propose to use Gaussian quadrature of order $q=4 \cdot p +1$, i.e.\ using~$n_q=2 \cdot p$ abscissae per interval $T \in \cT_h$. 

	\subsection{On Local Minimizers}
	\label{sec:local}
	Above, we proved that the global NLP minimizer converges to a (or the, in case it is unique) global  minimizer of \refOCP. However, practical NLP solvers can often only compute critical points, which have a local minimality certificate at best.
	For collocation methods, all critical  points of \refOCP have convergent critical NLP points if the mesh is sufficiently fine. For PBF the above global convergence result implies a more favorable assertion: For every strict local  minimizer of \refOCP there is exactly one convergent strict local NLP minimizer if the mesh is sufficiently fine. Below we explain the reason why this follows from the above global convergence property.
	
	Consider a strict local minimizer $\tilde{x}^\star$ of \refOCP. By definition of a local minimizer, inactive box constraints $x_L\leq x\leq x_R$ could be imposed on \refOCP such that $\tilde{x}^\star$ is the unique global minimizer of a modified problem. Upon discretization we would keep the box constraints as $\bx_L \leq \bx \leq \bx_R$. From the above convergence result, since $\tilde{x}^\star$ is unique with inactive box constraints, $\bx$ must converge to $\tilde{x}^\star$, leaving the NLP box constraints inactive, as if they had been omitted as in the original problem.
	
	\section{Solving the Nonlinear Program}\label{sec:solvingNLP}
	The core of our work is the above convergence proof for our Penalty-Barrier Finite Element (PBF) Method. Below, we point out one particular way in which the resulting NLP can be solved if the functions are sufficiently differentiable. This is compared to the Legendre-Gauss-Radau (LGR) collocation method  in terms of numerical cost and conditioning aspects.
	
	\subsection{Formulation of the NLP}
	We identify $x_h \in \cX_{h,p}$ with a finite-dimensional vector $\bx \in \R^\nX$ via nodal values (in particular, left LGR points), as depicted in Figure 1 with $\times$. Additional nodes $+$ for~$z_h$ could be used to allow for discontinuities, but to facilitate  easier comparisons between PBF and LGR this was not done here. We describe how the NLP can be solved numerically for $\bx$.
	
	Using a Gauss-Legendre quadrature rule of abscissae $s_j \in \overline{\Omega}$ and quadrature weights $\alpha_j \in \R_+$ for $j=1,\ldots,\nquad$ for PBF, we define two functions $F_h,C_h$ as 
		\begin{align}
		F_h(\bx) &:= \sum_{j=1}^{\nquad} \alpha_j \cdot f\big(\dot{y}_h(s_j),{y}_h(s_j),{z}_h(s_j),s_j\big),\label{eqn:NLP:F}\\
		C_h(\bx) &:= \begin{pmatrix}
		b\big(y_h(t_1),y_h(t_2),\dots,y_h(t_M)\big)\\[8pt]
		\sqrt{\alpha_1} \cdot c\big( \dot{y}_h(s_1),y_h(s_1),z_h(s_1),s_1\big)\\
		\sqrt{\alpha_2} \cdot c\big( \dot{y}_h(s_2),y_h(s_2),z_h(s_2),s_2\big)\\
		\vdots\\
		\sqrt{\alpha_\nquad} \cdot c\big( \dot{y}_h(s_\nquad),y_h(s_\nquad),z_h(s_\nquad),s_\nquad\big)
		\end{pmatrix}. \label{eqn:NLP:C}
		\end{align}
		Note that $ 	r_h(x_h) = \|C_h(\bx)\|_2^2. 	$
		The optimization problem in $\bx$ for PBF is 
		\begin{align}
		\operatornamewithlimits{min}_{\bx\in\mathbb{R}^{\nX}} \phi(\bx):= F_h(\bx) + \frac{1}{2 
			\omega} 
		\|C_h(\bx)\|_2^2 + \tau 
		\Gamma_h(\bx),\label{eqn:NLP}
		\end{align}
		where we emphasize that, due to an oversampling with quadrature points that exceed the dimension of $\bx$, there is no solution to $C_h(\bx)=0$, in general, hence why it is crucial to consider the penalty formulation \eqref{eqn:NLP} as is.
	
	In comparison, for LGR the NLP functions $F_h,C_h$ are exactly the same as for PBF, however with $n_q=p$ and $s_j,\alpha_j$ the LGR quadrature points and weights, and optionally a different weighting for the rows of $C_h$. Here, the optimization problem in $\bx$ for direct collocation  is
		\newcommand{\cI}{\mathcal{I}}
		\begin{align}
		\begin{split}
		\operatornamewithlimits{min}_{\bx\in\mathbb{R}^{\nX}} \qquad  &   F_h(\bx) \\
		\text{s.t.}                               \qquad &   C_h(\bx)=0\,,\quad Q_h(\bx)\geq \bO,
		\end{split}\label{eqn:NLPcol}
		\end{align}
		with $Q_h$ a function that evaluates $z_h$ at  collocation points.
		Note that, for collocation methods, it would not be possible in general	for the number $n_q$ of probing points for the path constraints
		to exceed the polynomial degree~$p$ of $\cX_{h,p}$, since otherwise there are too many constraints for too few degrees of freedom, resulting in an over-determined NLP.
	
	\subsection{Primal-Dual Penalty-Barrier NLP Solvers}
	Though \eqref{eqn:NLP} and \eqref{eqn:NLPcol} look different, they are actually similar in that both can be solved efficiently with the same type of numerical optimizations methods, that is, \textit{penalty-barrier} methods \cite{ForsgrenGill,ALM_IPM,SUMT,ForsgrenSIAMREV,NumOpt,BoydConvexOpt}.
	
	Penalty-barrier methods find numerical solutions that yield a KKT-residual on the order of magnitude of a user-specified tolerance value $\tol>0$. They work by reforming problems in the form \eqref{eqn:NLPcol} into a problem in the form \eqref{eqn:NLP} by devising a suitable logarithmic barrier function for the inequality constraints $Q$, and using $\omega=\tau=\tol$. For PBF, we therefore also propose choosing $\omega=\tau$ equal to $\tol$.
	
	When minimizing functions of the form $\phi$, the Hessian becomes ill-conditioned when $\omega,\tau$ are small, which is the practically relevant case.
		The mainstay of the practical efficiency of penalty-barrier methods is therefore their capability of circumventing the conditioning issues related to small penalty-barrier parameters. This is achieved in two ways:
		(i)~As is known for path-following methods~\cite{ForsgrenSIAMREV}, by minimizing $\phi$ successively for iteratively decreasing parameters $\omega,\tau$, the overall method's iteration count depends little on the final values of $\omega,\tau$.
		(ii) By using a \textit{primal-dual} formulation of the optimality conditions of~\eqref{eqn:NLP}, the optimality conditions match a dual regularized version of the optimality conditions of \eqref{eqn:NLPcol}, hence the condition number of the linearized regularized KKT systems is in fact better than for the unregularized KKT system. It is well-understood that under certain contraint qualifications the condition numbers in the KKT matrices are bounded and independent of $\omega,\tau$ \cite{ForsgrenGill,GMO14}.
	
	The method due to Forsgren and Gill \cite{ForsgrenGill} is one particular primal-dual penalty-barrier method, globalized via line search. 
		The method is widely considered as a simple, yet efficient prototype for modern large-scale NLP solvers. In their paper, the method is proposed for solving a problem in the form \eqref{eqn:NLPcol}. However, the method effectively minimizes the primal merit function given in~\cite[Eqn~1.1]{ForsgrenGill}, which matches exactly with our function $\phi$ in \eqref{eqn:NLP}. They review in detail the benign conditioning aspects of the linear system matrices involved, and describe how the iteration count becomes insensitive to $\omega,\tau$ (which, in their notation, are both replaced with a single number $\mu>0$).
	
	\subsection{Computational Cost}
	For the Forsgren-Gill method, we compare the cost between solving \eqref{eqn:NLP} from PBF and \eqref{eqn:NLPcol} from LGR.
	
	The Forsgren-Gill method is an iterative method, each iteration consisting of two steps: (i)~solving the primal-dual linear system, and (ii)~performing the line search. Assuming the number of iterations for PBF and LGR is similar (which is reasonable, since this rather depends on geometric behaviours of the underlying optimal control problem, as verified below in the numerical experiments, where we report the iteration counts), we only need to compare the cost per iteration.
	
	Since the line search is of negligible cost, only the costs for solving the linear system need to be compared. Due to their sparse structure, for optimal control problems it is attractive to form and solve the primal Schur complement, as proposed for a different context in \cite[Eqn~13]{ALM_IPM}:
		$$  \bS := \nabla_{\bx\bx}^2\cL(\bx,\blambda) + \frac{1}{\omega}\cdot\nabla_\bx C_h(\bx)\cdot \nabla_\bx C_h(\bx)\t,$$
		where $\blambda$ is the vector of Lagrange multipliers.
	
	Since the linear system is only in the primal variables after reduction, this matrix has the same dimension $\nX$ for~\eqref{eqn:NLP} and \eqref{eqn:NLPcol}. Further, the matrix has the same sparsity pattern for PBF and LGR. Namely, it is block-tridiagonal of small dense square overlapping blocks that have the dimension $n_x \cdot (p+1)$. For both methods, the number of blocks on the diagonal equals the total number of finite elements, thereby amounting to the dimension~$\nX$.
	
	The computational cost of evaluation and the number of non-zeros of the Hessian and Jacobian for PBF and LGR can be compared from the formula of the Lagrangian function $\cL(\bx,\blambda):=F_h(\bx)-\blambda\t \cdot C_h(\bx)$:
		\begin{multline}  
		\cL(\bx,\blambda) = \sum_{j=1}^{\nquad} \alpha_j \cdot \Bigg( f\big(\dot{y}_h(s_j),{y}_h(s_j),{z}_h(s_j),s_j\big) \\
		+ \blambda_j\t \cdot c\big(\dot{y}_h(s_j),{y}_h(s_j),{z}_h(s_j),s_j\big)\Bigg)
	    + \blambda_b\t \cdot b\big(y_h(t_1),y_h(t_2),\dots,y_h(t_M)\big) .
	    	\label{eqn:LagrangianFunction}
	    	\end{multline}
		$Q_h$ and $\Gamma_h$ have been omitted. This is because $Q_h$ is linear, with trivial Jacobian, namely a submatrix of an identity matrix. For interior-point methods, the barrier-function of $Q_h$ is the negative sum of the logarithms of its elements. $\Gamma_h$ in turn is just a weighted version of this sum, hence its Jacobian is a row-scaled version of the Jacobian of $Q_h$.
	
	Both LGR and PBF use the same formula for $\cL$, however PBF uses twice as many sample points $n_q$ for the constraints as LGR (as we propose).
		We acknowledge that $\blambda$ has no context in \eqref{eqn:NLP} other than by the definition of the Lagrangian, whereas  $\blambda$  constitutes the dual solution for \eqref{eqn:NLPcol}. For \eqref{eqn:NLP} Forsgren-Gill computes these to $\blambda=-\frac{1}{\omega}\cdot C_h(\bx)$.

	In summary, the computational cost  of obtaining all function values and the assembly of all derivative matrices, including $\bS$, is at most twice the computational cost for PBF compared to LGR. Accordingly, the number of nonzeros in any derivative matrices is  at most twice that for PBF as for LGR. In constrast, the cost for solving the linear system is identical for PBF and LGR when using the primal reduced matrix $\bS$.
	
	In the experiments below, we always report the bandwidth and number of nonzeros of $\bS$ for PBF and LGR, as well as the number of NLP iterations for the Forsgren-Gill method to converge.

	\subsection{Derivatives and Sparsity}
	Second-order NLP solvers need $\nabla_\bx\cL(\bx,\by)$, $\nabla^2_{\bx,\bx}\cL(\bx,\by)$ and $\nabla_\bx C(\bx)$.
	We see that \eqref{eqn:NLP:F}--\eqref{eqn:NLP:C} essentially use the same formulas as collocation methods, except that they replace collocation points with quadrature points $\tau_j$, and weight the constraints $C$ by $\sqrt{\alpha_j}$. Hence, the Jacobian of PBF has the same sparsity structure as for a collocation method of same degree $p$, except the Jacobian in PBF has more rows than in LGR; e.g.\ twice as many rows as a collocation method of the same degree if we choose to use Gauss-Legendre quadrature with $q=2p$ abscissae per element, as proposed.
	The observations for $\nabla\t_\bx C(\bx)$ can be extended to $\nabla^2_{\bx,\bx}\cL(\bx,\by)\in \R^{n \times n}$: PBF($p$) and LGR($p$) of the same degree $p$ on the same mesh have the same Hessian dimension and sparsity structure. Section~\ref{sec:numExp:Sparsity} presents sparsity plots for PBF(5) and LGR(5).
	
\newcommand{\abbFig}{Fig.}
\newcommand{\abbfig}{Fig.}
	
	\section{Numerical Experiments}
	\label{sec:numerical}
	The scope of this paper is the transcription method.
	Practical aspects in solving the NLP \eqref{eqn:PBP_h} are discussed in \cite{CDC20PBF}, where we also show non-zero patterns of the sparse Jacobian and Hessian of the constraints and Lagrangian function; and show that the computational cost roughly compares to solving the NLPs from LGR collocation. Below, we present numerical results for two test problems when using our transcription method and minimizing \eqref{eqn:PBP_h} for the given instance, mesh, and finite element degree.
	
	
	\subsection{Convex Quadratic Problem with a Singular Arc}\label{sec:numExp:CQP}
\begin{figure}[tb]
		\centering
		\includegraphics[width=0.6\columnwidth]{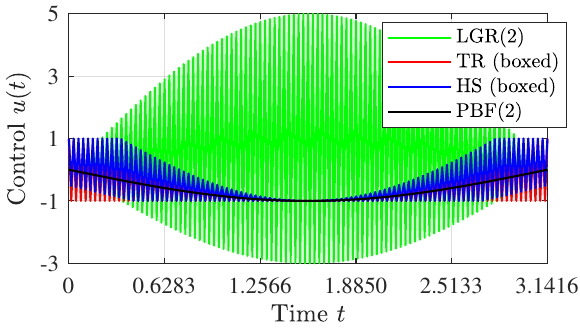}
		\caption{Comparison of control solutions for three collocation methods and PBF.}
		\label{fig:CQVP}
	\end{figure}

	Consider a small test problem, which demonstrates convergence of PBF in a case where direct collocation methods ring:
	\begin{equation}
		\begin{aligned}
			&\min_{y,u} \quad &&\int_0^{\pi} \left(\,y_0(t)^2 + \cos^{(1-m)}(t) \cdot  u(t)\,\right)\,\mathrm{d}t,\\
			&\text{s.t.} &&\dot{y}_{k-1}(t)=y_{k}(t),\text{ for }k=1,\dots,m\,,\qquad\dot{y}_m(t)=u(t),
		\end{aligned}
	\end{equation}
	where $\cos^{(1-m)}$ is the $(1-m)^{th}$ derivative of $\cos$, with negative derivative meaning antiderivative. Figure~\ref{fig:CQVP} shows the numerical solutions for $m=1$ for Trapezoidal (TR), Hermite-Simpson (HS), LGR collocation and PBF, where the latter two use polynomial degree $p=2$. For higher degree $p$, LGR would still ring when $m\geq p-1$. TR and HS require box constraints on $u$ for boundedness of their NLP minimizers.
	
	\subsection{Second-Order Singular Regulator}
	This bang-singular control problem from \cite{SOSC} with \mbox{$t_E=5,$}\,\mbox{$\eta=1$} is given as
	\begin{equation}
	\begin{aligned}
	&\min_{y,u} &\quad &\int_0^{t_E} \left(\,y_2(t)^2 + \eta \cdot  y_1(t)^2\,\right)\,\mathrm{d}t,\\
	&\text{s.t.} & y_1(0)&=0,\quad y_2(0)=1,\\
	&& \dot{y}_1(t)&=y_2(t),\quad \dot{y}_2(t)=u(t),\\
	&& -1\leq u(t)&\leq 1.
	\end{aligned}\label{eqn:SOSR_ocp}
	\end{equation}
	
	Both LGR and PBF use $100$ elements of degree $p=5$. $\bS$~has 24507 non-zeros and bandwidth 15 for both discretizations. Forsgren-Gill solves PBF in 40 and LGR in 41 NLP iterations.
	
\begin{figure}[tb]
	\centering
	\includegraphics[width=0.6\columnwidth]{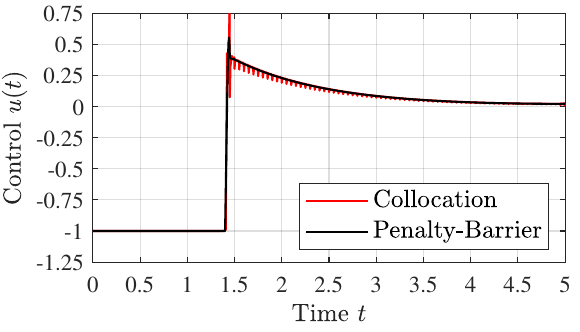}
	\caption{Comparison of control solutions from LGR and PBF for the Second-Order Singular Regulator.}
	\label{fig:SOSController}
	\end{figure}

	Figure~\ref{fig:SOSController} presents the control profiles of the two numerical solutions. LGR shows ringing on the time interval $[1.5,\,5]$ of the singular arc. In contrast, PBF converges with the error $\|u^\star(t)-u_h(t)\|_{L^2([1.5,\,5])} \approx 1.5 \cdot 10^{-4}$.
	
	\subsection{Aly-Chan Problem}
	The problem in \cite{AlyChan}, namely \eqref{eqn:SOSR_ocp} with \mbox{$t_E=\pi/2,$} \mbox{$\eta=-1,$} has a smooth totally singular control.
	
\begin{figure}[tb]
		\centering
		\includegraphics[width=0.6\columnwidth]{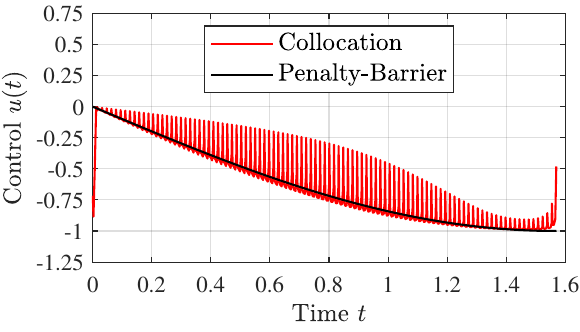}
		\caption{Comparison of control solutions from LGR and PBF for the Aly-Chan Problem.}
		\label{fig:NumExp_AlyChan}
	\end{figure}

	Both LGR and PBF use $100$ elements of degree $p=5$. $\bS$ has the same nonzero pattern as before. Forsgren-Gill for PBF/LGR converges in 48/43 iterations.
	Figure~\ref{fig:NumExp_AlyChan} presents the control profiles of the two numerical solutions. PBF converges, with error
	$  \|u^\star(t)-u_h(t)\|_{L^2(\Omega)} \approx 3.7 \cdot 10^{-6}     . $
	LGR does not converge for this problem; cf.~\cite[Fig.~3]{Kameswaran},\cite{ChenBiegler16}.
	
	\subsection{Regular state-constrained problem}\label{sec:numExp:Sparsity}
\begin{figure}
		\centering
		\includegraphics[width=0.6\columnwidth]{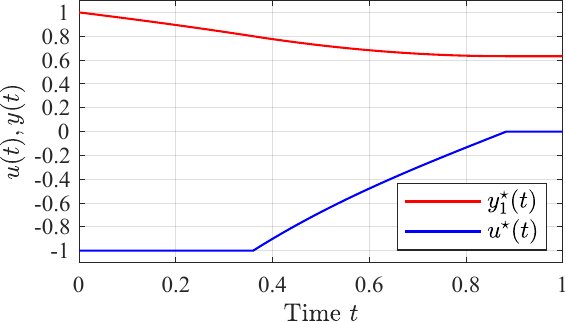}
		\caption{Analytical solution to \eqref{eqn:ExampleOCP2}.}
		\label{fig:exp2}
	\end{figure}

	We now consider a test problem for which both types of methods converge with success, so that we can compare conditioning, convergence, and rate of convergence to a known analytical solution:
	\begin{align*}
		\min_{y,u}& & J&= y_2(1),\label{eqn:ExampleOCP2}\\
		\text{s.t.}& & y_1(0)&=1,\ \, \dot{y}_1(t)=\frac{u(t)}{2y_1(t)},\ \, \sqrt{0.4}\leq y_1(t),\\
		& & y_2(0)&=0,\ \, \dot{y}_2(t)=4 y_1(t)^4 + u(t)^2,\ \,  -1\leq u(t)\,.
	\end{align*}
	The solution is shown in \abbfig~\ref{fig:exp2}.
	$u^\star$ is constant outside $t_0=1-\frac{\sqrt{41}}{10}\approx 0.35$ and $t_1=t_0+\log 2-\frac{\log(\sqrt{41}-5)}{2}\approx 0.88$, between which $u^\star(t)= 0.8 \sinh\left(2 (t - t_1)\right)$, yielding $J\approx 2.0578660621682771255864272367598$. 
	
\begin{figure}
		\centering
		\includegraphics[width=0.6\columnwidth]{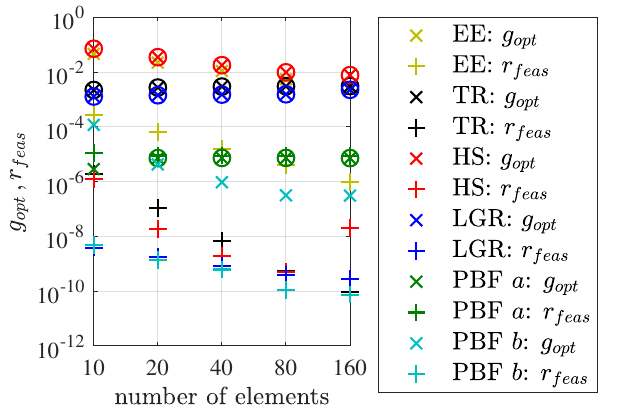}
		\caption{Convergence of optimality gap and feasibility residual
			LGR and PBF use polynomial degree $p=5$. PBF uses two different values: a) $\omega=10^{-5}$, b) $\omega=10^{-10}$.}
		\label{fig:loglogxpl2}
	\end{figure}

	All methods yield accurate solutions. \abbFig~\ref{fig:loglogxpl2} shows the convergence of the optimality gap and feasibility residual of a respective method. Remarking on the former, we computed $J^\star-J(x_h)$ and encircled the cross when $J(x_h)<J^\star$. Note in the figure that for $\geq 40$ elements the most accurate solutions in terms of feasibility are found by PBF with $\omega=10^{-10}$. Further, we find that the collocation methods significantly underestimate the optimality value for this experiment.
	
	Now we discuss rates of convergence. Convergence of only first order is expected because higher derivatives of $y^\star$ are non-smooth and $u^\star$ has edges. Indeed, $r_\text{feas}$ converges linearly for all methods. PBF5 with $\omega=10^{-5}$ stagnates early because it converges to the optimal penalty solution, which for this instance is converged from $20$ elements onwards. $g_\text{opt},r_\text{feas}$ are then fully determined by $\omega$. The issue is resolved by choosing $\omega$ smaller. LGR5 and PBF5 with $\omega=10^{-10}$ converge similarly, and stagnate at $r_\text{feas}\approx 10^{-10}$. Due to the high exponent in the objective, small feasibility errors in the collocation methods amount to significant underestimation of the objective.
	
\begin{figure}
		\centering
		\includegraphics[width=0.6\columnwidth]{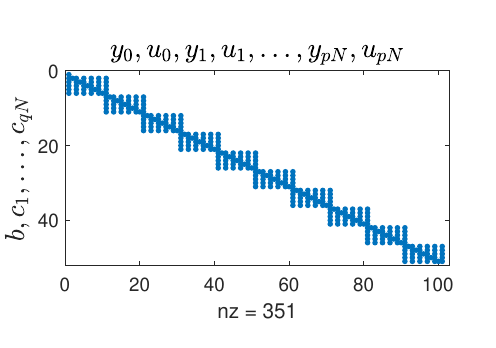}
		\caption{Sparsity of $\nabla_\bx C(\bx)\t$ for LGR5 when $h=\frac{1}{10}$, i.e.~$N=10$. For LGR, notice $q=p-1$. The discretization does not depend on $u_{pN}$.}
		\label{fig:sparsitylgr}
	\end{figure}

\begin{figure}
		\centering
		\includegraphics[width=0.6\columnwidth]{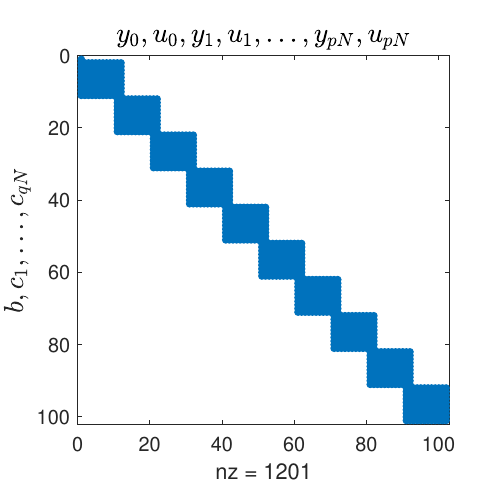}
		\caption{Sparsity of $\nabla_\bx C(\bx)\t$ for PBF5 when $h=\frac{1}{10}$, i.e.~$N=10$, with $q=2p$.}
		\label{fig:sparsitypbf}
	\end{figure}

	Finally, we look into computational cost. Solving the collocation methods with IPOPT and the PBF5 discretization with the interior-point method in \cite{ForsgrenGill}, the optimization converges in $\approx 20$ iterations for any discretization. Differences in computational cost can arise when one discretization results in much denser or larger problems than others. Here, we compare the sparsity structure of the Jacobian $\nabla_\bx C(\bx)\t$ for LGR5 in \abbfig~\ref{fig:sparsitylgr} and PBF5 in \abbfig~\ref{fig:sparsitypbf}, each using a mesh size of $h=\frac{1}{10}$.
	Note that for PBF5, $C(\bx)$ has more rows in the Jacobian than LGR5, thus the Jacobian has hence more non-zeros. However, critical for computations is the primal Schur complement $\mathbf{\Sigma}=\nabla_{\bx\bx}^2 \cL(\bx,\blambda)+\nabla_\bx C(\bx)\t \bD \nabla_\bx C(\bx)$, which is used when solving the KKT system via the reduced form, where $\bD$ is a diagonal matrix. $\mathbf{\Sigma}$ is a narrow-banded matrix with dense band of the same bandwidth for LGR5 and PBF5.
	
	With regard to computational cost, it follows from Fig.~\ref{fig:loglogxpl2} that the ability to choose $\omega$ in PBF can be advantageous. In particular, on coarse meshes, one may opt for small feasibility residual by manually decreasing $\omega$, whereas with a collocation method one is stuck with the feasibility residual that one obtains for that particular mesh. The figure shows this: For $\omega=10^{-10}$, even on the coarsest mesh the PBF method achieves a solution that has a smaller feasibility residual than other methods on the same mesh. For this problem this becomes possible because the path constraint could be satisfied with zero error by choosing $y$ a polynomial of degree 3 (because here PBF uses $p=5$).

	\subsection{Van der Pol Controller}
	
	This problem uses a controller to stabilize the van der Pol differential equations on a finite-time horizon. The problem is stated as
	\begin{equation*}
	\begin{aligned}
	&\min_{y,u} &\quad \frac{1}{2}\cdot &\int_0^4 \left(\,y_1(t)^2 + y_2(t)^2\,\right)\,\mathrm{d}t,\\
	&\text{s.t.} & y_1(0)&=0,\quad\quad y_2(0)=1,\\
	&& \dot{y}_1(t)&=y_2(t),\\
	&& \dot{y}_2(t)&=-y_1(t)+y_2(t) \cdot \left(\,1-y_1(t)^2\,\right) + u(t),\\
	&& -1 & \leq u(t)\leq 1. 
	\end{aligned}
	\end{equation*}
	The problem features a bang-bang control with a singular arc on one sub-interval. The discontinuities in the optimal control are to five digits at $t_1 = 1.3667$ and $t_2 = 2.4601$.
	
	We solved this problem with LGR collocation on $100$ uniform elements of order $5$. We compare this solution to the one obtained with PBF using $100$ uniform elements of order $p=5$, with $\omega=10^{-10}$ and $\tau=10^{-10}$.
	
\begin{figure}[tb]
		\centering
		\includegraphics[width=0.6\columnwidth]{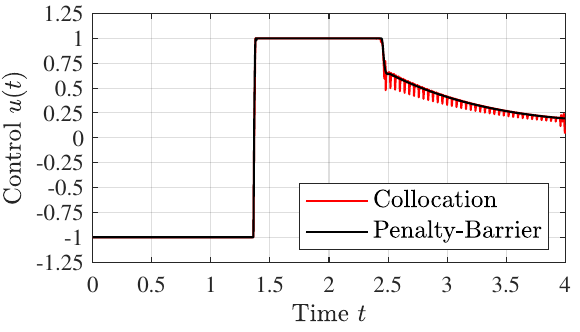}
		\caption{Comparison of control input obtained with Penalty-Barrier method against LGR collocation  for the Van der Pol problem.}
		\label{fig:NumExp_VanderPol}
	\end{figure}

	Figure~\ref{fig:NumExp_VanderPol} presents the control profiles of the two numerical solutions.
	LGR shows ringing on the time interval $[t_2,\,4]$ of the singular arc. In contrast,  PBF converges to the analytic solution. The solution satisfies the error bounds
	$e(0)\approx 7.0\cdot 10^{-2}$, $e(t_2)\approx 1.2\cdot 10^{-2}$, $e(2.5)\approx 8.17\cdot 10^{-4}$, and $e(2.6)\approx 9.6 \cdot 10^{-5}$, where
	$e(\hat{t}):=\|u^\star(t)-u_h(t)\|_{L^2([\hat{t},4])}$. The larger errors in the vicinity of the jumps occur due to the non-adaptive mesh.
	
	\subsection{Reorientation of an Asymmetric Body}
	This nonlinear problem from  \cite[Ex.\ 6.12, eqn.\ 6.123]{Betts2nd} in the ODE formulation is considered numerically challenging for its minimum time objective, the control appearing only linearly in the dynamics, the problem having a bang-bang solution and multiple local minima with identical cost. Since the solution is bang-bang, rates of convergence (at least for the optimality gap) can only be linear. We compare convergence of three collocation methods and PBF, where a polynomial degree $p=3$ is used for the $hp$-methods.
	
	Using the same initial guess from forward integration of an approximate control solution, LGR and PBF 
	converge on average in 200 iterations. LGR was solved with IPOPT and PBF was solved with a penalty-barrier interior-point method presented in \cite{ForsgrenGill}. Both NLP solvers cost one linear system solve per iteration. For $\omega=\tau=10^{-3}$, the finite element solution $x^\star_h$ converges to the penalty-barrier minimizer $x^\star_{\omega,\tau}$ sooner, which however is not very feasible for the DOP at hand. The other collocation methods' NLPs were also solved using IPOPT, which terminated on local infeasibility for TR and HS. In contrast, LGR and PBF provide  numerical solutions that converge at similar rates, which stagnate around $10^{-6}$ for the feasibility residual and $10^{-4}$ for the optimality gap. These methods converge at similar rates. 
	The size of the differential constraint violation and optimality gap for HS, LGR and PBF for different mesh sizes are given  in Table~\ref{tab:r_feas}. For HS, due to box constraints on end-time, which has been expressed as the first state, the optimality gap is negative and equal to the lower box constraint on the first state.
	
	As with any regularization method, good values for $\omega,\tau$ can be found on the fly by saving each barrier solution and terminating when differences between subsequent barrier solutions stop decreasing. For computation of the gap, we determined $J^\star:=28.6298010321$ from PBF(3) on 2048 elements, where $r_\textrm{feas}\approx 2.4e-7$.
	
	\newcommand{\linebreakcell}[2][c]{\begin{tabular}[#1]{@{}c@{}}#2\end{tabular}}
	\newcommand{\lbqcell}[2]{\linebreakcell{#1/\\ \hspace{2mm}#2}\xspace}
	
	\begin{table}[tb]
	\caption{$L^2(\Omega)$-norm for differential constraints / optimality gap of the Asymmetric Body Reorientation Problem. All methods use consistency order $p=3$. $N_\text{el}$ is the number of elements.}
		\label{tab:r_feas}
		\centering
        \begin{tabular}{|c||c|c|c|c|c|c|c|}
\hline
$N_\text{el}$ 	& HS 						& LGR 						& \linebreakcell{PBF\\$\omega=\tau=10^{-3}$} 	& \linebreakcell{PBF\\ $\omega=\tau=10^{-7}$} 	& \linebreakcell{PBF\\ $\omega=\tau=10^{-10}$} 	\\ \hline\hline
8   			& \lbqcell{3.2e-2}{-1.3e-1} & \lbqcell{6.1e-3}{1.2e+0} 	& \lbqcell{4.5e-3}{-4.4e-1} 				& \lbqcell{1.9e-4}{1.4e+0}  				& \lbqcell{1.8e-4}{1.4e+0} 					\\ \hline
32  			& \lbqcell{3.6e-3}{-1.3e-1} & \lbqcell{6.6e-5}{2.2e-2} 	& \lbqcell{4.2e-3}{-4.5e-1}	 				& \lbqcell{2.5e-5}{3.6e-2}  				& \lbqcell{7.7e-6}{4.2e-1} 					\\ \hline
128 			& \lbqcell{1.4e-3}{-1.3e-1} & \lbqcell{1.0e-6}{5.7e-4} 	& \lbqcell{4.1e-3}{-4.5e-1} 				& \lbqcell{1.7e-6}{6.9e-4}  				& \lbqcell{3.4e-7}{9.5e-3} 					\\ \hline
512 			& \lbqcell{7.0e-4}{-1.3e-1} & \lbqcell{2.1e-6}{9.2e-5} 	& \lbqcell{4.1e-3}{-4.5e-1} 				& \lbqcell{1.4e-6}{7.4e-5}  				& \lbqcell{1.1e-8}{6.1e-4} 					\\ \hline
        	\hline
        \end{tabular}
	
	\end{table}

	\subsection{\submithla{Obstacle Avoidance Problem}}
	Since we limited our presentation to a convergence analysis for global minimizers, we give this example to demonstrate PBF's practical capability to also converge to non-global minimizers.
	
\begin{figure}[tb]
		\centering
		\includegraphics[width=0.6\columnwidth]{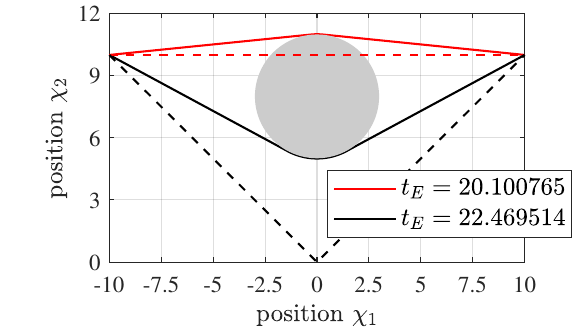}
		\caption{Optimal trajectories from PBF for the Obstacle Avoidance Problem.}
		\label{fig:NumExp_Obst}
	\end{figure}

	Consider the minimum-time trajectory from $\vec{\chi}_0=[-10\ 10]\t$ to $\vec{\chi}_E=[10\ 10]\t$ around an obstacle at $\vec{\chi}_C=[0\ 8]\t$ of radius $R=3$:
		\begin{equation*}
		\begin{aligned}
		\min_{\vec{\chi},u,t_E} \quad t_E&,\\
		\text{s.t.} \quad \vec{\chi}(0)&=\vec{\chi}_0,\quad\vec{\chi}(t_E)=\vec{\chi}_E, \quad\|\vec{\chi}(t)-\vec{\chi}_C\|_2^2\geq R^2\\
		\dot{\vec{\chi}}(t)&=\big[\cos\big(u(t)\big)\ \sin\big(u(t)\big)\big]\t,
		\end{aligned}
		\end{equation*}
		Passing the obstacle above or below results in two locally optimal trajectories. Both are found by PBF, depicted in Figure~\ref{fig:NumExp_Obst}, using the dashed curves as initial guesses (with $t_E$ and $u$ computed feasible from $\vec{\chi}$ via integration and differentiation, respectively) on $100$ finite elements of degree 5.
	
	The computed times as in the figure are accurate except to the last digit. The red/black trajectory converges in 52/58 NLP iterations. For comparison, LGR of the same degree and mesh converges in 51/51 iterations. $\bS$ has 73521 nonzeros and bandwidth 25 for both PBF and LGR.
	
	\subsection{Pendulum in Differential-Algebraic Form}
	In this example from~\cite[Chap.~55]{betts2015collection}, a control force decelerates a frictionless pendulum to rest. The objective is to minimize the integral of the square of the control:
	\begin{equation*}
	\begin{aligned}
	&\min_{\vec{\chi},\xi,u} & &\int_0^{3} u(t)^2\,\mathrm{d}t,\\
	&\text{s.t.} & \vec{\chi}(0)&=[1\ 0]\t,\quad \dot{\vec{\chi}}(0)=\vec{0},\\
	&&\vec{\chi}(3)&=[0\ -1]\t,\quad \dot{\vec{\chi}}(3)=\vec{0},\\
	&& \ddot{\vec{\chi}}(t)&=[0\ -9.81]\t +2 \cdot \vec{\chi}(t) \cdot \xi(t) + \vec{\chi}^\perp(t) \cdot u(t),
	\end{aligned}
	\end{equation*}
	with an additional DAE constraint introduced below. The ODE for $\vec{\chi}$ is a force balance in the pendulum mass. $u(t)$ is the control force acting in the direction $\vec{\chi}^\perp := [-\chi_2\ \chi_1]\t$.
	
	The DAE constraint determines the beam force $\xi(t)$ in the pendulum arm in an implicit way, such that the length remains 1 for all time;
	\cite[Chap.~55]{betts2015collection} uses
	\begin{align}
	0&=\|\dot{\vec{\chi}}(t)\|_2^2 -2 \cdot \xi(t) -g \cdot \chi_2(t).\label{eqn:BeamDAE1}
	\end{align}
	The following alternative constraint achieves the same:
	\begin{align}
	\tag{\ref{eqn:BeamDAE1}'}
	0&=\|\vec{\chi}(t)\|_2^2 - 1\,.\label{eqn:BeamDAE3}
	\end{align}
	\eqref{eqn:BeamDAE1} is a DAE of index 1, whereas \eqref{eqn:BeamDAE3} is of index~3. 

	In the following we study the convergence of TR, HS, LGR ($p=5$) and PBF ($p=5$) on meshes of increasing size. Here, the collocation methods are solved with IPOPT in ICLOCS2, whereas PBF is solved with Forsgren-Gill as before.
	
	TR \& HS are likely to converge at a slower rate than PBF \& LGR. However, our focus is primarily on determining whether a given method converges, and only secondarily on rates of convergence. To find out where solvers struggle, we consider three variants of the pendulum problem,
	\begin{enumerate}[\text{\ Case\,}A]
		\item where we consider the original problem with \eqref{eqn:BeamDAE1} as given in \cite{betts2015collection}.
		\item where we add the path constraint $\xi(t)\leq 8$.
		\item where we exchange \eqref{eqn:BeamDAE1} with \eqref{eqn:BeamDAE3}.
	\end{enumerate}
	
\begin{figure}[tb]
		\centering
		\includegraphics[width=0.6\columnwidth]{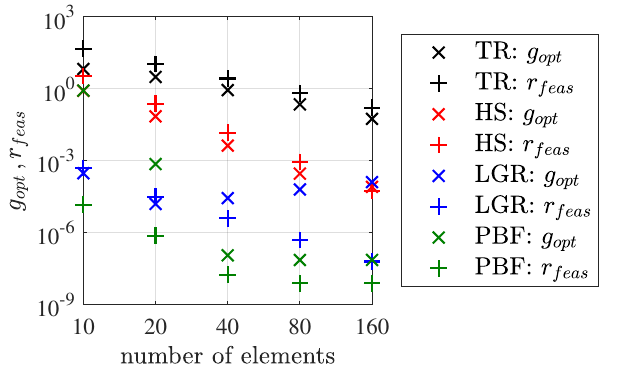}
		\caption{Convergence of optimality gap and feasibility residual for Pendulum example, case A.}
		\label{fig:pendulumcaseAlog}
	\end{figure}

	All methods converge  for case A. Figure~\ref{fig:pendulumcaseAlog} shows that TR converges slowly, while HS, LGR and PBF converge fast. At small magnitudes of $g_{\text{opt}},r_{\text{feas}}$, further decrease of LGR and PBF deteriorates, presumably due to limits in solving the NLP accurately under rounding errors.
	
\begin{figure}[tb]
		\centering
		\includegraphics[width=0.6\columnwidth]{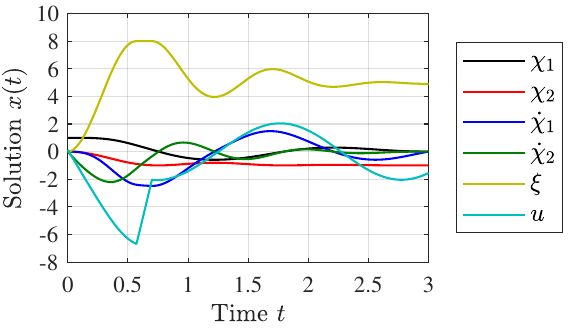}
		\caption{Numerical solution of PBF on $80$ elements for Pendulum example, case B.}
		\label{fig:pendulumcaseBarcs}
	\end{figure}

\begin{figure}[tb]
		\centering
		\includegraphics[width=0.6\columnwidth]{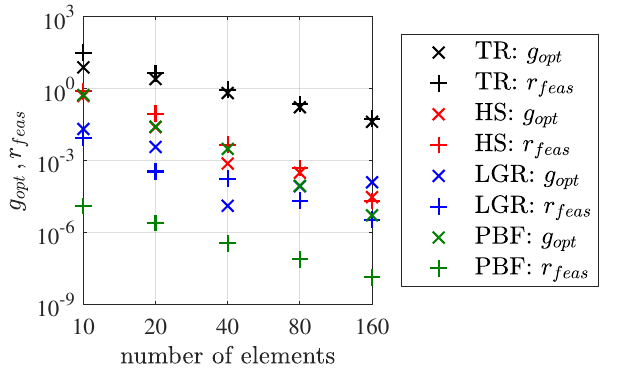}
		\caption{Convergence of optimality gap and feasibility residual for Pendulum experiment, case B.}
		\label{fig:pendulumcaseBlog}
	\end{figure}

	Case B is shown in Figure~\ref{fig:pendulumcaseBarcs}. The control force decelerates the pendulum more aggressively before the pendulum mass surpasses the lowest point, such that the beam force obeys the imposed upper bound. Figure~\ref{fig:pendulumcaseBlog}  confirms convergence for all methods. The rate of convergence is  slower compared to case A, as expected, because  the solution of $u$ is locally non-smooth.

	
\begin{figure}[tb]
		\centering
		\includegraphics[width=0.6\columnwidth]{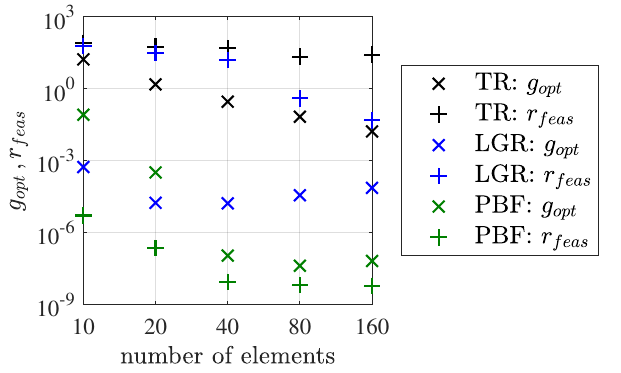}
		\caption{Convergence of optimality gap and feasibility residual for Pendulum example, case C.}
		\label{fig:pendulumcaseClog}
	\end{figure}

\begin{figure}[tb]
		\centering
		\includegraphics[width=0.6\columnwidth]{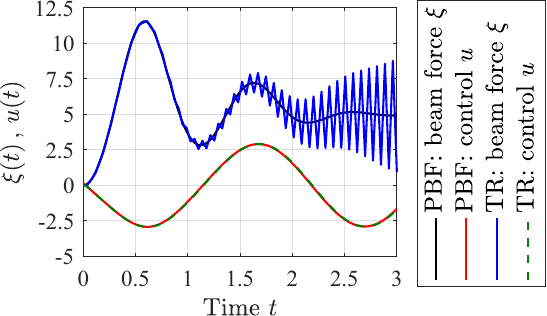}
		\caption{Numerical solutions of PBF and TR on $80$ elements for Pendulum example, case C. The optimal control is identical to case A.}
		\label{fig:pendulumcaseCarcs}
	\end{figure}

	For case~C, some collocation methods struggle: For HS on all meshes, the restoration phase in IPOPT converged to an infeasible point, indicating infeasibility of \eqref{eqn:NLPcol} for this scheme \cite[Sec.~3.3]{IPOPT}. For TR, the feasibility residual does not converge, as shown in Figure~\ref{fig:pendulumcaseClog}.
	Figure~\ref{fig:pendulumcaseCarcs} shows that this is due to ringing in the numerical solution for the beam force.
	Regarding LGR, Figure~\ref{fig:pendulumcaseClog} shows that the feasibility residual converges only for relatively fine meshes. In contrast to the collocation methods, PBF converges as fast as for case~A.
	
	Finally, we discuss the computational cost: Using 160 elements, the PBF discretization results in $\bS$ of bandwidth~30, with 44964 nonzeros for cases~\mbox{A,\,C} and 42951 nonzeros for case~B; requiring 66, 51, and 66 NLP iterations for cases A--C. LGR yields the same sparsity pattern for $\bS$ as PBF, solving on average in $30$ IPOPT iterations (with second-order corrections).
		
	\section{Conclusions}
	\label{sec:conclusions}
	We presented PBF and proved convergence under mild and easily-enforced assumptions. Key to the convergence proof is the formulation of a suitable unconstrained penalty-barrier problem, which is discretized using finite elements and solved with primal-dual penalty-barrier NLP solvers.
	
	Theorem~\ref{thm:order} provides high-order convergence guarantees even if the component $z$ has discontinuities, provided that the trajectory can be approximated accurately in the finite element space; see~\eqref{eqn:InfBound} and the discontinuous elements in Figure~\ref{fig:finiteelementsyz}. It is  a practical matter to employ an adaptive meshing technique for achieving this in an efficient manner.
	
	While this paper has a theoretical focus, the practicality of our novel transcription has been illustrated in numerical examples. The scheme converged for challenging problems, which included solutions with singular arcs and discontinuities. These problems caused issues for  three commonly used direct transcription methods based on collocation, namely TR, HS and LGR.  

	\section*{Acknowledgements}
	The authors thank Yuanbo Nie for help and advice on numerical experiments in ICLOCS2.
	
	\appendix
	
	\section{Converting a Bolza Problem into~\refOCP}
	\label{sec:convert}
	
	Many control, state or parameter estimation problems can be converted into the
	Bolza form~\cite{Betts2nd,Rawlingsetal:2017} below, possibly over a variable time domain and/or with a terminal cost function~$f_E$ (also known as a Mayer term or final cost):
	\begin{equation*}
		\operatornamewithlimits{min}_{\chi,\upsilon,\xi,\tau_0,\tau_E} \int_{\tau_0}^{\tau_E} \, f_r\left(\dot{\chi}(\tau),\chi(\tau),\upsilon(\tau),\xi,\tau\right)\mathrm{d}\tau + f_E(\chi(\tau_E),\tau_E)
	\end{equation*}
	subject to
	\begin{align*}
		b_B\left(\chi(\tau_0),\chi(\tau_E\right),\tau_0,\tau_E)&= 0,\\
		c_e\left(\dot{\chi}(\tau),\chi(\tau),\upsilon(\tau),\xi,\tau\right)&=0 \text{ f.a.e.\ }\tau \in (\tau_0,\tau_E), \\
		c_i\left(\dot{\chi}(\tau),\chi(\tau),\upsilon(\tau),\xi,\tau_0,\tau_E,\tau\right) &\leq 0 \text{ f.a.e.\ }\tau \in (\tau_0,\tau_E),
	\end{align*}
	where the state is $\chi$, the weak derivative  $\dot{\chi}:=\mathrm{d}\chi/\mathrm{d}\tau$ and   the input is $\upsilon$. Note that the starting point $\tau_0$, end point $\tau_E$ and a constant vector of parameters $\xi\in\R^{n_\xi}$ are included as optimization variables.  The above problem can be converted into the Lagrange form~\refOCP as follows.
	
	Move the Mayer term into the  integrand by noting that $f_E(\chi(\tau_E),\tau_E) = \phi(\tau_E)= \int_{\tau_0}^{\tau_E}\dot{\phi}(\tau)\mathrm{d}\tau$
	if we let
	$
	\phi(\tau):=f_E(\chi(\tau),\tau)
	$
	with initial condition $\phi(\tau_0)=0$.  Recall that for minimum-time problems, we usually let $f_E(\chi(\tau),\tau):=\tau$ and $f_r:=0$, so that $\dot{\phi}(\tau)=1$.
	
	By introducing the auxiliary  function $s$, convert the inequality constraints into the equality constraints
	$
	s(\tau)+c_i\left(\dot{\chi}(\tau),\chi(\tau),\upsilon(\tau),\xi,\tau_0,\tau_E,\tau\right) = 0
	$
	and inequality constraint $s\geq 0$. Introduce the auxiliary functions $(\upsilon^+,\upsilon^-)\geq 0$ with the substitution $\upsilon = \upsilon^+ - \upsilon^-$ so that  the algebraic variable function is defined as $z:=(s,\upsilon^+,\upsilon^-)$.
	
	The problem above with a variable domain is converted onto a fixed domain with  $t\in(0,1)$ via the transformation $\tau = \tau_0+(\tau_E-\tau_0) \cdot t$ so that $t_0:=0$, $t_E=1$.
	
	Define the  state for problem~\refOCP as $y:=(\chi,\phi,\xi,\tau_0,\tau_E)$ and introduce additional equality constraints in order to force ($\mathrm{d}\xi/\mathrm{d}t,\,\mathrm{d}\tau_0/\mathrm{d}t,\,\mathrm{d}\tau_E/\mathrm{d}t)=0$.
	
	The expressions for $f,\, c,\, b$ can now be derived using the above.

	\section{Lebesgue Equivalence for Polynomials}\label{sec:Appendix_LebesgueIdentity}
	Let $T:=(a,b) \in \cT_h$ and $p \in \N_0$. We show that
	$$ 	\|\beta\cdot u\|_{L^\infty(T)} \leq \frac{p+1}{\sqrt{|T|}} \cdot \|\beta \cdot u\|_{L^2(T)} \quad \forall u \in \cP_p(T),\forall \beta \in \R. $$
	
\begin{figure}[tb]
		\centering
		\includegraphics[width=0.6\columnwidth]{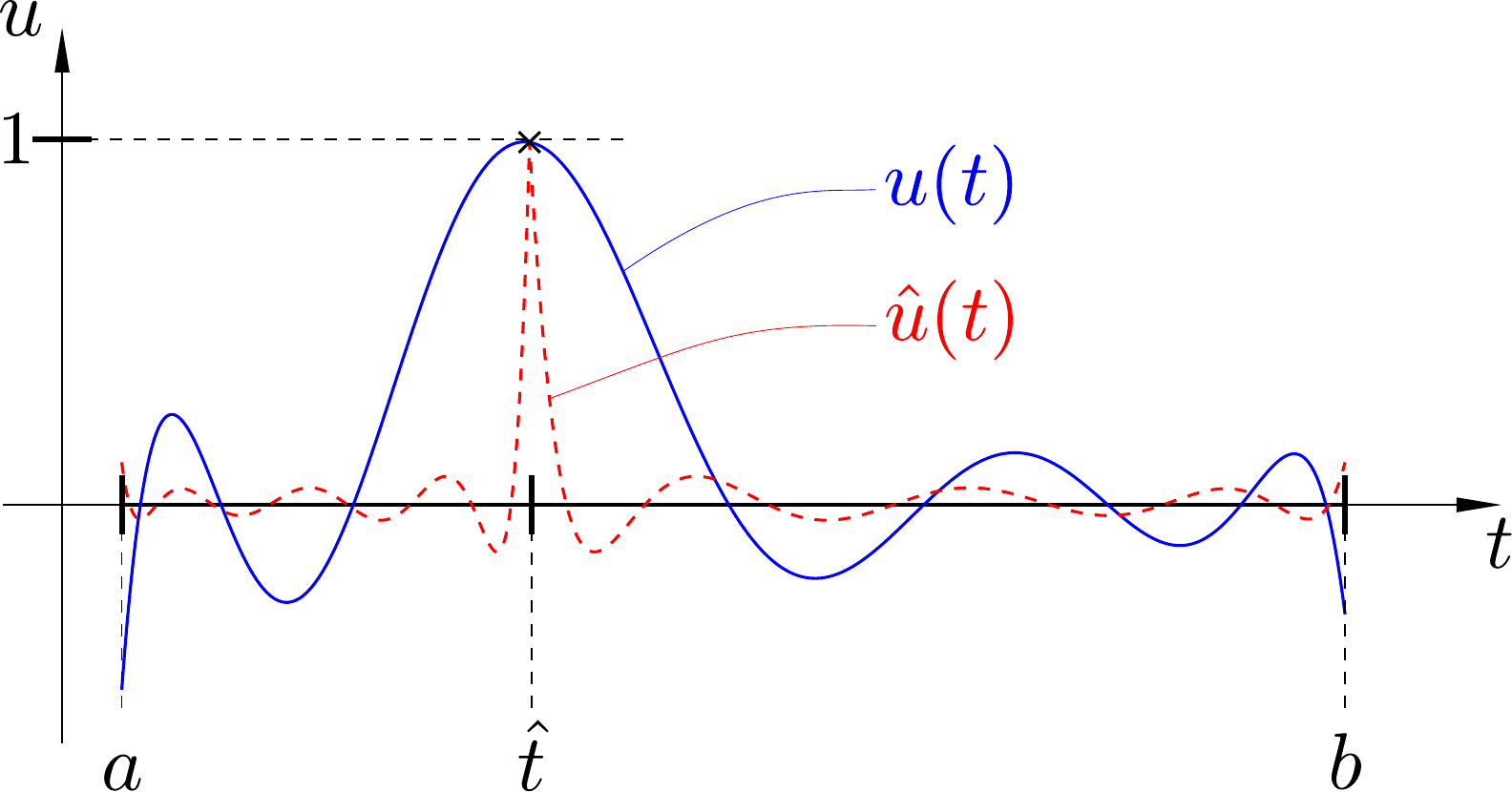}
		\caption{Polynomial $u$ and piecewise polynomial $\hat{u}$ over $T$.}
		\label{fig:polyplothur}
	\end{figure}

	Choose $u \in \cP_p(T)$ arbitrary. Since $\|\beta\cdot u \|_{L^k(T)} = |\beta| \cdot \|u\|_{L^k(T)}$ holds for both $k \in \lbrace 2,\infty\rbrace$, and for all $\beta\in\R$, w.l.o.g.\ let $\|u\|_{L^\infty(T)}=1$.
	Since $\operatorname{sgn}(\beta)$ is arbitrary, w.l.o.g.\ let $u(\hat{t})=1$ for some $\hat{t} \in \overline{T}$.
	Define $T_L:=[a,\hat{t}]$, $T_R:=[\hat{t},b]$, $\hat{\cP}_p:=\cP_p(T_L)\cap\cP_p(T_R)\cap\cC^0(T)$,
	and
	$ 	\hat{u} := \operatornamewithlimits{arg\,min}_{v \in \hat{\cP}_p}\big\lbrace \|v\|_{L^2(T)}\ \big\vert\ v(\hat{t})=1 \big\rbrace. 	$
	Since $\cP_p \subset \hat{\cP}_p$, it holds $\|u\|_{L^2(T)}\geq \|\hat{u}\|_{L^2(T)}$.
	Figure~\ref{fig:polyplothur} illustrates $u,\hat{u}$ for $p=8$.
	
	\newcommand{\hur}{{\hat{u}_\text{ref}}\xspace}
	\newcommand{\Tr}{{T_\text{ref}}\xspace}
	Use $\|\hat{u}\|^2_{L^2(T)}=\int_a^b \hat{u}(t)^2\mathrm{d}t=(b-a)/2\cdot\int_{-1}^1 \hur(t)^2\mathrm{d}t=\frac{|T|}{2} \cdot \|\hur\|^2_{L^2(\Tr)}$, where $\hur$ is $\hat{u}$ linearly transformed from $T$ onto $\Tr:=(-1,1)$.
	Since $\|\hat{u}\|_{L^2(T)}$ is invariant under changes of $\hat{t}$ because $\hat{u}(\hat{t}+(b-\hat{t})\cdot\xi)=\hat{u}(\hat{t}+(\hat{t}-a)\cdot\xi)$ $\forall \xi \in [0,1]$, w.l.o.g.\ we can assume for $\hat{u}$ that $\hat{t}=b$ and hence $\hur(1)=1$.
	Since minimizing the $L^2(\Tr)$-norm, $\hur$ solves
	\begin{equation}
	\label{eqn:CQP_hur}
	\operatornamewithlimits{min}_{v \in \cP_p(\Tr)}\ 1/2 \cdot \int_\Tr v(t)^2\,\mathrm{d}t\text{ subject to } v(1)=1.
	\end{equation}
	
	We represent $\hur=\sum_{j=0}^p \alpha_j \cdot \phi_j$, where $\phi_j$ is the $j^\text{th}$ Legendre polynomial. These satisfy \cite{refLegendre}: $\phi_j(1)=1\quad\forall j\in \N_0,\quad \int_\Tr \phi_j(t)\cdot\phi_k(t)\mathrm{d}t=\delta_{j,k} \cdot \gamma_j\quad\forall j,k\in\N_0,$
	where $\gamma_j:=2/(2\cdot j +1)$ and $\delta_{j,k}$ the Kronecker delta. We write $\bx=(\alpha_0,\alpha_1,\dots,\alpha_p)\t\in\R^{p+1}$, ${D}=\opdiag(\gamma_0,\gamma_1,\dots,\gamma_p)\in\R^{(p+1)\times (p+1)}$ and $\be\in\R^{p+1}$.
	Then \eqref{eqn:CQP_hur} can be written in~$\bx$:
	$$ 	\min_{\bx\in\R^{p+1}}\ \psi(\bx):=1/2\cdot\bx\t\cdot{D}\cdot\bx \text{ subject to }\be\t\cdot\bx=1. 	$$
	From the optimality conditions \cite[p.~451]{NumOpt}\away{
		$
		\left[\begin{array}{c|c}
		{D} & \be\\
		\hline
		\phantom{{}^{{}^A}}\be\t\phantom{{}^{{}^A}} & 0
		\end{array}\right]
		\cdot
		\left(\begin{array}{c}
		\bx\\
		\hline
		-\lambda
		\end{array}\right) =
		\left(\begin{array}{c}
		\bO\\
		\hline
		1
		\end{array}\right)$} follows $\bx={D}\inv\cdot\be\cdot\lambda$ and $\be\t\cdot{D}\cdot\be\cdot\lambda=1$.
	Using $\be\t\cdot{D}\cdot\be=\sum_{j=0}^p \gamma_j=\frac{(p+1)^2}{2}$ yields $\lambda=1/(p+1)^2$ and
	$\psi(\bx)=\frac{1}{2} \cdot ({D}\inv\cdot\be\cdot\lambda)\t\cdot{D}\cdot({D}\inv\cdot\be\cdot\lambda)=\frac{\lambda}{2}=\frac{1}{(p+1)^2}$.
	Hence, $\frac{1}{2}\cdot\|\hur\|^2_{L^2(\Tr)}=1/(p+1)^2$.
	Hence, $\frac{1}{2}\cdot\|\hat{u}\|^2_{L^2(T)}=\frac{|T|}{2}\cdot 1/(p+1)^2$.
	Hence, $\|u\|_{L^2(T)}\geq\|\hat{u}\|_{L^2(T)}=\sqrt{|T|}/(p+1)\cdot \underbrace{\|u\|_{L^\infty(T)}}_{=1}$,
	or, \mbox{$\|\beta \cdot u\|_{L^2(T)}\geq \sqrt{|T|}/(p+1)\cdot \|\underbrace{\beta \cdot u}_{\tilde{u}}\|_{L^\infty(T)}$}.
	In conclusion:
	\begin{align}
		\|\tilde{u}\|_{L^\infty(T)}\leq \frac{p+1}{\sqrt{|T|}}\cdot\|\tilde{u}\|_{L^2(T)}\ \, \forall \tilde{u} \in \cP_p(T)\, \forall T \in \cT_h. \label{eqn:PropAppendix2}
	\end{align}
	
	\section{Order of Approximation for Non-smooth and Continuous Non-differentiable Functions}\label{app:3}
	
	In the following we illustrate that the assumption $\ell>0$ in \eqref{eqn:InfBound} is rather mild. To this end, we consider two pathological functions for $g:=x^\star_{\omega,\tau}$. In our setting, $n_y=0,\,n_z=1$, and we  interpolate a given pathological function $g$ with $x_h \in \cX_{h,p}$ over $\Omega=(-1,1)$. We use $p=0$.
	
	\paragraph{A function with infinitely many discontinuities}
	The first example is a non-smooth function that has infinitely many discontinuities. Similar functions can arise as optimal control solutions; cf.\ Fuller's problem \cite{Fuller}.
	
\begin{figure}[tb]
		\centering
		\includegraphics[width=0.6\columnwidth]{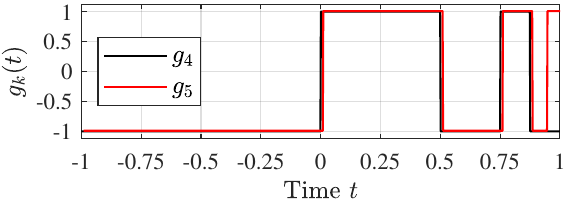}
		\caption{Nested step-function $g_k$ for $k=4,\,5$\,.}
		\label{fig:NestedStep}
	\end{figure}

	Consider the limit $g_\infty$ of the following series:
	\newcommand{\sign}{\operatorname{sign}\xspace}
	\begin{align*}
		g_0(t) := -1\,,\qquad	g_{k+1}(t) := \left\lbrace \begin{matrix}
			g_{k}(t) & \text{if }t \leq 1-2^{-k}\\
			-g_k(t) & \text{otherwise}
		\end{matrix} \right.& &k=0,1,2,\dots.
	\end{align*}
	$g_\infty$ switches between $-1$ and $1$ whenever~$t$ halves its distance to $1$. Figure~\ref{fig:NestedStep} shows $g_k$ for $k=4,\,5$\,.
	
	Using mesh-size $h=2^{-k}$ for some $k \in \N$,  define $u(t):=g_k(t)\in \cX_{h,p}$. Hence,
	$	\operatornamewithlimits{inf}_{x_h \in \cX_{h,p}}
	\|g_\infty-x_h\|_\cX 
	\leq \|g_\infty-u\|_{L^2(\Omega)} $.
	It follows that
	$$ 	|u(t) - g_\infty(t)| \leq \left\lbrace \begin{matrix}
	0 & \text{if }t \leq 1-2^{-k}\\
	2 & \text{otherwise}
	\end{matrix} \right. $$
	Hence, $ \|g_\infty-u\|_{L^2(\Omega)}\leq \|g_\infty-u\|_{L^1(\Omega)} \leq 2/2^k = \cO(h^1).$
	Therefore, all $\ell \in (0,0.5]$ satisfy~\eqref{eqn:InfBound}.
	
	\paragraph{A continuous but nowhere differentiable function}
	Consider the following Weierstrass function, which is continuous but non-differentiable:
	$$ 	g(t) := \frac{1}{2} \cdot \sum_{k=0}^\infty a^{k} \cdot \cos(7^k \cdot \pi \cdot t)$$
	for $0<a\leq 0.5$. This function with range $\subset[-1,1]$ satisfies the H\"older property
	$ |g(t)-g(s)| \leq C \cdot |t-s|^\alpha $
	with some $C \in \R_+$ for $\alpha = -\log(a)/\log(7)$\, \cite{Zygmund}.
	For $a\leq 0.375$ we have $\alpha\geq 0.504$\,.
	
	According to this property, a piecewise constant interpolation $u \in \cX_{h,p}$ of $g$ satisfies
	$ |g(t)-u(t)| \leq |g(t)-g(s)| \leq C \cdot |t-s|^\alpha\leq |h|^\alpha $. 
	In conclusion, $\operatornamewithlimits{inf}_{x_h \in \cX_{h,p}}\left\lbrace \|g-x_h\|_\cX \right\rbrace \leq \|g-u\|_{L^2(\Omega)} \leq \|g-u\|_{L^1(\Omega)} = \cO(h^\alpha)$. Therefore, all $\ell \in (0,\alpha-0.5]$ satisfy~\eqref{eqn:InfBound}.
	
	\section{Proof of Lemma~\ref{lem:BoundLipschitz_Fr}}\label{sec:Appendix_ProofLemma1}
	
		The boundedness follows from (A.2).
		
		Lipschitz continuity of  $r$ is not as straightforward. We will make use of the following trace theorem \cite{TraceDing}: For an open interval $I \subseteq \Omega$ it holds that $\|u\|_{L^2(\partial I)} \leq K \cdot \|u\|_{H^1(I)}$	with a constant $K$ independent of $u$. Assume $|u|$ attains its essential supremum on $\overline{\Omega}$ at $t=t^\star$. Choosing $I=(t^\star,t_E)\subset\Omega$, then
		$  \|u\|_{L^\infty(\Omega)} = |u(t^\star)| \leq \|u\|_{L^2(\partial I)}.  $
		Using this together with the above bound and
		$  \|u\|_{H^1(I)}\leq\|u\|_{H^1(\Omega)}  $
		results in
		\begin{align}
			\|u\|_{L^\infty(\Omega)} \leq K \cdot \|u\|_{H^1(\Omega)}.  \label{eqn:TraceResult}
		\end{align}
		
		Below, for a generic Lipschitz continuous function $g:\R^k\rightarrow \R^{n_g}$ with Lipschitz-constant $L_g$ and $\|\cdot\|_1$-bound $|g|_\text{max}$, we use the relation
		\begin{align*}
		    \begin{split}
			&\left|\|g(\xi_2)\|_2^2\,-\,\|g(\xi_1)\|_2^2\right|
			=\big|\|g(\xi_2)\|_2+\|g(\xi_1)\|_2\big| \cdot \underbrace{\big|\|g(\xi_2)\|_2-\|g(\xi_1)\|_2\big|}_{\leq\|g(\xi_2)-g(\xi_1)\|_2}\\
			&\leq n_g \cdot \big|\|g(\xi_2)\|_1+\|g(\xi_1)\|_1\big| \cdot \|g(\xi_2)-g(\xi_1)\|_1
			\leq n_g \cdot 2 \cdot |g|_{\text{max}} \cdot L_g \cdot \|\xi_2 - \xi_1\|_1\,,
		    \end{split}\label{eqn:LipQuadPenaltyFunc}
		\end{align*}
		where we used $|\alpha^2-\beta^2|=|\alpha+\beta|\cdot|\alpha-\beta|$ in the first line and the triangular inequality in the second line.
		\away{
			\begin{align*}	
				&\|g(\xi_2)-g(\xi_1)\|_2 \cdot \|g(\xi_2)-g(\xi_1)\|_2\\
				&\quad \leq \dim(g) \cdot \|g(\xi_2)-g(\xi_1)\|_1 \cdot \|g(\xi_2)-g(\xi_1)\|_1\\
				&\quad \leq \dim(g) \cdot L_g \cdot \|\xi_2-\xi_1\|_1 \cdot (\|g(\xi_2)\|_1+\|g(\xi_1)\|_1)
			\end{align*}
		}
		Using the above bound, we can show Lipschitz continuity of~$r$:
		\begin{align*}
			&|r(x_2)-r(x_1)|\leq \int_\Omega \Big| \left\|c\left(\dot{y}_2(t),{y}_2(t),z_2(t),t\right)\right\|_2^2 - \left\|c\left(\dot{y}_1(t),{y}_1(t),z_1(t),t\right)\right\|_2^2\Big| \mathrm{d}t\\
			&\quad + \Big|\left\|b\left(y_2(t_1),\ldots,y_2(t_M)\right)\right\|_2^2
			-\left\|b\left(y_1(t_1),\ldots,y_1(t_M)\right)\right\|_2^2\Big| \\
			&\leq \int_\Omega 2 \cdot \nc \cdot |c|_\text{max} \cdot L_c \cdot \left\|\begin{pmatrix}
				\dot{y}_2(t)-\dot{y}_1(t)\\
				y_2(t)-y_1(t)\\
				z_2(t)-z_1(t)
			\end{pmatrix}\right\|_1 \mathrm{d}t\\
			&\quad + 2 \cdot \nb \cdot |b|_\text{max} \cdot L_b \cdot \underbrace{\left\|\begin{pmatrix}
					y_2(t_1)-y_1(t_1)\\
					\vdots\\
					y_2(t_M)-y_1(t_M)
				\end{pmatrix}\right\|_1}_{\leq M \cdot \|y_2-y_1\|_{L^\infty(\Omega)}} \\
			&\leq 2  \cdot \nc \cdot |c|_\text{max} \cdot L_c \cdot \left(\|\dot{y}_2-\dot{y}_1\|_{L^1(\Omega)}+\|y_2-y_1\|_{L^1(\Omega)}+\|z_2-z_1\|_{L^1(\Omega)}\right)\\
			&\quad + 2  \cdot \nb \cdot |b|_\text{max} \cdot L_b \cdot M \cdot K \cdot \|y_2-y_1\|_{H^1(\Omega)},
		\end{align*}
		where \eqref{eqn:TraceResult} has been used 
		to bound $\|y_2-y_1\|_{L^\infty(\Omega)}$.
		
		If $y_2=y_1$ then we see the result shows Lipschitz continuity of $r$ with respect to $\|z\|_{L^1(\Omega)}$. Using
		$$  \|u\|_{L^1(\Omega)} \leq \sqrt{|\Omega|} \cdot \|u\|_{L^2(\Omega)}\quad \forall u\in L^1(\Omega)$$
		according to \cite[Thm.~2.8, eqn.~8]{Adams}, and the definition of $\|\cdot\|_\cX$, we arrive at
		\begin{align*}
			&\|\dot{y}_2-\dot{y}_1\|_{L^1(\Omega)} +
			\|y_2-y_1\|_{L^1(\Omega)} +
			\|z_2-z_1\|_{L^1(\Omega)} \\
			\leq& \sqrt{|\Omega|}\cdot\left(\|\dot{y}_2-\dot{y}_1\|_{L^2(\Omega)} +
			\|y_2-y_1\|_{L^2(\Omega)} +
			\|z_2-z_1\|_{L^2(\Omega)}\right)
			\leq 3 \cdot \sqrt{|\Omega|}\cdot\|x_2-x_1\|_\cX,
		\end{align*}
		which shows Lipschitz continuity of $r$ with respect to~$\|x\|_\cX$.
		
		Lipschitz continuity of $F$ follows from Lipschitz continuity of $f$:
		\begin{align*}
			|F(x_2)-F(x_1)|
			&\leq \int_\Omega |f(\dot{y}_2(t),y_2(t),z_2(t))-f(\dot{y}_1(t),y_1(t),z_1(t))|\,\mathrm{d}t\\
			&\leq \int_\Omega L_f \cdot \left\|\begin{pmatrix}
				\dot{y}_2(t)-\dot{y}_1(t)\\
				{y}_2(t)-{y}_1(t)\\
				{z}_2(t)-{z}_1(t)
			\end{pmatrix}\right\|_{1}\mathrm{d}t
			\leq L_f \cdot \left\|\begin{pmatrix}
				\dot{y}_2-\dot{y}_1\\
				{y}_2-{y}_1\\
				{z}_2-{z}_1
			\end{pmatrix}\right\|_{L^1(\Omega)}
		\end{align*}
	
	\section{Properties of the $\log$-Barrier Function}\label{sec:Appendix_BarrierFunctioProperties}
	Let $0<\zeta\ll 1$ be a fixed small arbitrary number.
	\begin{lem}[Order of the $\log$ Term]\label{lem:tauLw}It holds:
		$ 	\left|\tau \cdot \log\left({\tau}/{L_\omega}\right)\right| = \cO\left(\tau^{1-\zeta}\right).	$
	\end{lem}
	\begin{proof}
		We use $L_\omega$ from \eqref{eqn:LipLw}, where $L_F\geq 2,\ L_r\geq 2$ and $0<\tau\leq\omega\leq 1$.  We get
		\begin{align*}
			\left|\tau\cdot\log\left({\tau}/{L_\omega}\right)\right|
			&= \tau \cdot \left(|\log(\tau)-\log(L_\omega)|\right) \leq \tau \cdot \left(|\log(\tau)|+|\log(L_\omega)|\right)\\
			&= \tau \cdot \left(\left|\log\left(L_F+\frac{L_r}{2\cdot\omega}\right)\right|+|\log(\tau)|\right)\\
			&\leq \tau \cdot \left( 1+|\log(L_F)|+\left|\log\left(\frac{L_r}{2\cdot\omega}\right)\right|+|\log(\tau)| \right)\\
			&\leq \tau \cdot\Big( \underbrace{1+|\log(L_F)| +|\log(L_r/2)|}_{= \cO(1)}+\underbrace{|\log(\omega)|}_{\leq|\log(\tau)|}+|\log(\tau)|\Big)\\
			&= \cO(\tau) + \cO(\tau \cdot |\log(\tau)|).
		\end{align*}
		In the third line above, we used the fact that for $\alpha,\beta\geq 2$,
		follows $\log(\alpha+\beta)\leq\log(\alpha)+\log(\beta).$
		The result follows from $\tau \cdot |\log(\tau)| = \cO(\tau^{1-\zeta})$ by L'H\^opital:
		\begin{small}
			\begin{align*}
				\lim\limits_{\tau \rightarrow 0}\frac{\tau \cdot \log(\tau)}{\tau^{1-\zeta}}
				=\lim\limits_{\tau \rightarrow 0} \frac{\log(\tau)}{\tau^{-\zeta}}
				\myeqLH\lim\limits_{\tau \rightarrow 0}\frac{\frac{1}{\tau}}{-\zeta \cdot \tau^{-\zeta-1}}
				= \lim\limits_{\tau \rightarrow 0} \frac{\tau^\zeta}{-\zeta}
				=0
			\end{align*}
		\end{small}
	\end{proof}
	
	\begin{lem}[Bound for $\Gamma$]\label{lem:Gamma}
		If $x \in \cX$ with $\|z\|_{L^\infty(\Omega)} = \cO(1)$, then
		\begin{align*}
			\left|\tau \cdot \Gamma\left(\bar{x}\right)\right| &= \cO\left(\tau^{1-\zeta}\right)\,,&
			\left|\tau \cdot \Gamma\left(\check{x}\right)\right| &= \cO\left(\tau^{1-\zeta}\right)\,.
		\end{align*}
	\end{lem}
	\begin{proof}
		Since the definitions are similar, we only show the proof for $\bar{x}$:
		\begin{align*}
|\tau \cdot \Gamma(\bar{x})| \leq\left|\tau \cdot \sum_{j=1}^{n_z}\int_\Omega\,\log\left(\bar{z}_{[j]}(t)\right)\mathrm{d}t \right| \leq n_z \cdot |\Omega| \cdot \operatornamewithlimits{max}_{1\leq j\leq n_z} \|\tau \cdot \log(\bar{z}_{[j]})\|_{L^\infty(\Omega)}\\
\leq n_z \cdot |\Omega| \cdot \Big(\underbrace{\cO\left(\tau^{1-\zeta}\right)}_{\text{bound for }\bar{z}_{[j]}< 1} + \underbrace{\cO(\tau)}_{\text{bound for }\bar{z}_{[j]}\geq 1}\Big)
= \cO\left(\tau^{1-\zeta}\right).
			\label{eqn:aux:PBP_OptGap}
		\end{align*}
		In the third line, we  distinguished two cases, namely
		$ 	\left|\log\left(\bar{z}_{[j]}(t)\right)\right| 	$
		attains its essential supremum at a $t\in \overline{\Omega}$ where either $\bar{z}_{[j]}(t)<1$ (case 1) or where $\bar{z}_{[j]}(t)\geq 1$ (case 2). In the first case, we can use  Lemma~\ref{lem:StrictInteriorness}~\&~\ref{lem:tauLw}. In the second case, we simply bound the logarithm using
		$	\|\bar{z}_{[j]}\|_{L^\infty(\Omega)}\leq\|z\|_{L^\infty(\Omega)}= \cO(1)$ to arrive at the term $\cO(\tau)$.
	\end{proof}

	
	
	\bibliographystyle{siamplain.bst}
	%
	
	\balance
	
	\bibliography{main}
	
	%
\end{document}